\documentclass[11pt]{article}

\usepackage[a4paper, left=2cm, right=2cm, top=2cm, bottom=2cm]{geometry} %space on sides
\usepackage{cmbright}
\usepackage[OT1]{fontenc}
\usepackage[english]{babel}
\usepackage{etex} %workspace
\usepackage{amsmath}
\usepackage{amsthm}
\usepackage{amsfonts}
\usepackage{amssymb}
\usepackage{wrapfig}
\usepackage{xcolor}  %for colored texts and individual colors 
\usepackage{extarrows}  %for larger and special arrows
\usepackage[pdfborder={0 0 0}]{hyperref}  %for Hyperlinks, pdfborder={0,0,0} removes red borderd hyperlinks 
\usepackage{enumitem} %enumeration with roman letters
\usepackage{graphicx} %to include graphics via "\includegraphics"
\usepackage{array} %for thiker hlines in tables (cf. newcommand below)
\usepackage[font=footnotesize,labelfont=bf]{caption} %for captions of graphics not included as "figure"
\usepackage{longtable}
\usepackage{pdfcomment} %for comments in PDF Output
\usepackage{caption}
\usepackage{subcaption}
\usepackage{changes} %for crossed out text
\usepackage{pdfpages} %including pdf documents as pages
\usepackage{datetime}
\usepackage{todonotes}
\usepackage{oplotsymbl} %in-line squares, pentagons & hexagons
\usepackage[backend=bibtex,style=numeric,maxbibnames=99]{biblatex}
\AtEveryBibitem{
	\clearfield{issn} %remove issn
	\clearfield{doi} %remove doi
	\ifentrytype{online}{}{ %remove url except for @online
		\clearfield{url}
	}
}
\usepackage{tikz-cd}
		\usetikzlibrary{calc} 
		\usetikzlibrary{intersections}
        \usetikzlibrary{arrows}
\usepackage{tabularx}
\usepackage{dsfont} %for \one-Symbol (indicator function)

\usepackage[bottom]{footmisc}

\newtheorem{thm}{Theorem}[section] % [section] sets theorem-counter within section
\newtheorem{lem}[thm]{Lemma}

\theoremstyle{definition} % for non-italic style

\theoremstyle{remark}
\newtheorem{remark}[thm]{Remark}

\newcommand{\EE}{\mathbb{E}}

\newcommand{\PP}{\mathbb{P}}

\newcommand{\var}{\operatorname{var}}

\newcommand{\dint}{\mathrm{d}}

\newcommand{\FindStat}[1]{
	\ifx&#1&
	\url{www.findstat.org}    % #1 is empty
	\else \url{www.findstat.org/#1}    % #1 is nonempty
	\fi}

\newcommand{\tc}{\operatorname{TC}}
\def\one{\mathds1}
\newcommand{\loweru}{\underline{u}}
\newcommand{\upperu}{\overline{u}}
\newcommand{\midarrow}{\tikz \draw[-triangle 90] (0,0) -- +(.1,0);}
\makeatletter
\let\@fnsymbol\@alph
\makeatother

\begin{document}

% ====================================================================
% TITLE & AUTHORS
% ====================================================================
	
\title{Vertex number of the typical cell in a\\ tri-directional Poisson line tessellation}

% author one information
\author{
Nils Heerten\footnote{Faculty of Mathematics, Ruhr University Bochum. Email: nils.heerten@rub.de}
\and
Janina H\"ubner\footnote{Faculty of Mathematics, Ruhr University Bochum \emph{and} Department of Economics and Sustainability in Health Care, Hochschule für Gesundheit. Email: janina.huebner@rub.de}
\and
Christoph Th\"ale\footnote{Faculty of Mathematics, Ruhr University Bochum. Email: christoph.thaele@rub.de}}

\date{}

\maketitle

\begin{abstract}
    This paper deals with the typical cell in a Poisson line tessellation in the plane whose directional distribution is concentrated on three equally spread values with possibly different weights. Such a random polygon can only be a triangle, a quadrilateral, a pentagon or a hexagon. The probability for each of these cases is determined explicitly in terms of the weights. Extremal cases are discussed as well.\\[2mm]
    \textbf{2020 Mathematics Subject Classification (MSC):} 60D05
 \end{abstract}

% ====================================================================
\section{Introduction and main result}
\label{sec:introduction_&_main_result}
% ====================================================================

\begin{figure}[bt]
	\centering
	\includegraphics[scale=0.2]{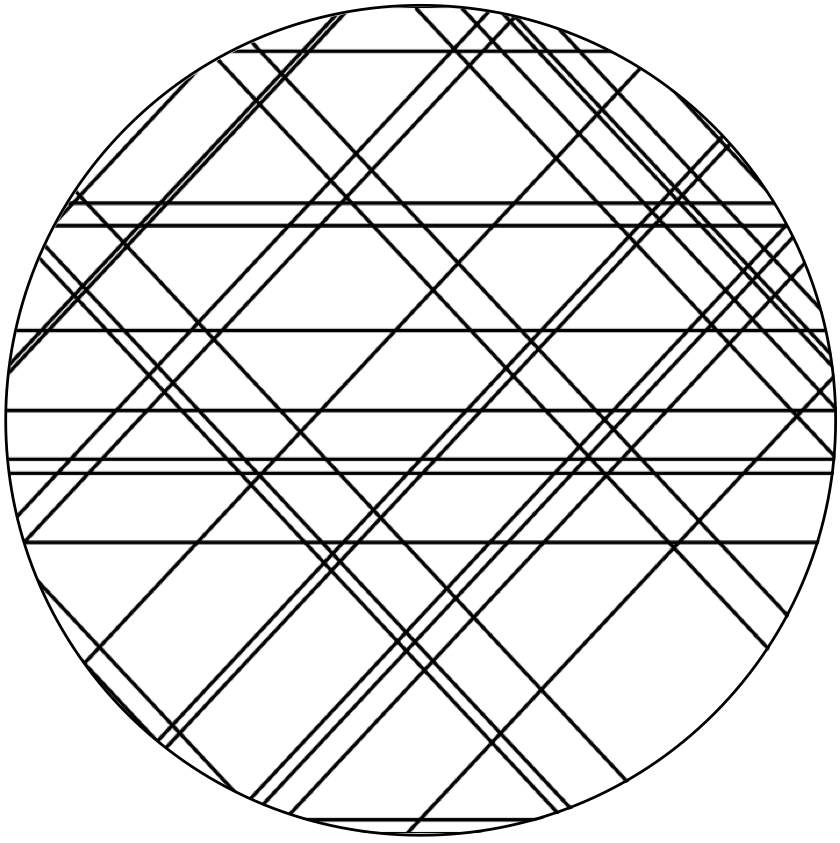}
	\caption[width=0.5\columnwidth]{A realization of a Poisson line tessellation with directional distribution $G_{1/3,1/3}$.}
	\label{fig:simulation_PLT_three_directions}
\end{figure}

Random tessellations are among the most classical objects to be studied in stochastic geometry. In this paper, we concentrate on random tessellations in the plane which are induced by stationary Poisson line processes. This infinite collection of random lines decomposes the plane into an infinite aggregate of random polygons. We recall from \cite[Chapter 9.5]{stoyan_kendall_mecke_2013} that the distribution of a Poisson line process is uniquely determined by an intensity parameter $\gamma\in(0,\infty)$ and a directional distribution $G$, which for us is a probability measure on $[0,\pi)$ satisfying $G(\theta)<1$ for all $\theta\in[0,\pi)$. We  are interested in the vertex number of the typical cell $\tc_G$ of a Poisson line tessellation with directional distribution $G$. Since this is an affine invariant quantity, we can and will from now on assume without loss of generality that the intensity satisfies $\gamma=1$. Informally, the typical cell can be thought of as a random polygon sampled `uniformly at random' from the infinite aggregate of all polygons induced by the Poisson line tessellation, regardless of its size an shape. Formally, 
\begin{equation}\label{eq:TypicalCell}
\PP(\tc_G\in\,\cdot\,) = \frac{1}{\EE\sum\limits_{C:m(C)\in[0,1]^2}1}\,\EE\sum_{C:m(C)\in[0,1]^2}\one_{\{C-m(C)\in\,\cdot\,\}},
\end{equation}
where we sum over all cells $C$ of the tessellation with the property that its lexicographically smallest vertex $m(C)$ is contained in the unit square $[0,1]^2$ (or any other Borel set with positive and finite Lebesgue measure). The systematic study of these random polygons dates back to works of Miles \cite{Miles64,miles_1973} in the isotropic case, where $G=G_{\rm unif}$ is the uniform distribution, and to that of George \cite{george_1987} and Mecke \cite{Mecke95} for general $G$. Our main focus lies on the vertex number $N_G$ of the typical cell $\tc_G$. Whereas it is well known that $\EE N_{G}=4$ for all directional distributions $G$, see \cite[Equation (9.70)]{stoyan_kendall_mecke_2013}, much less is known about the probabilities $\PP(N_{G}=n)$ for $n\in\{3,4,\ldots\}$, which in turn seem to depend on $G$ in a rather subtle way. In the isotropic case, $\PP(N_{G_{\rm unif}}=3)=2-{\pi^2\over 6}$ was determined by Miles \cite{Miles64} and $\PP(N_{G_{\rm unif}}=4)=\pi^2\log2-\frac13-\frac{7\pi^2}{36}-\frac72\sum_{j=1}^\infty\frac{1}{j^3} \approx 0.381466$ has been calculated by Tanner \cite{Tanner83}. For $n\geq 5$ there are only involved integral formulas and numerical results in \cite{Calka03} as well as the tail asymptotics from \cite{BonnetCalkaReitzner}, see also \cite[Table 9.2]{stoyan_kendall_mecke_2013}. For a class of discrete directional distributions $G$ the precise value for the triangle probability $\PP(N_G=3)$ has been determined in \cite{heerten_krecklenberg_thaele_2022}. The first non-trivial member in this class is the directional distribution $G_{p,q}$ concentrated on three equally spread angles. It is given by 
\begin{equation}
	G_{p,q}:=p\delta_{0}+q\delta_{\pi/3}+(1-p-q)\delta_{2\pi/3},
	\label{eq:g_3}
\end{equation}
where $\delta_{(\,\cdot\,)}$ denotes the Dirac-measure and the weights $0<p,q<1$ are such that $p+q<1$. A simulation of a Poisson line tessellation with directional distribution $G_{1/3,1/3}$ is shown in Figure \ref{fig:simulation_PLT_three_directions}.

Our main contribution is a complete description of the distribution of the vertex number $N_{p,q}:=N_{G_{p,q}}$ of the typical cell $\tc_{p,q}:=\tc_{G_{p,q}}$ in a Poisson line tessellation with directional distribution $G_{p,q}$. It should be observed that for this particular choice of $G$, the vertex number is a random variable concentrated on $\{3,4,5,6\}$. We emphasize that, as far as we are aware of, our result is the first complete distributional description of the vertex number of the typical cell in a Poisson line tessellation.

\begin{thm}
\label{thm:main_theorem}
    Let $N_{p,q}$ be the vertex number of the typical cell $\tc_{p,q}$ of a Poisson line tessellation with directional distribution $G_{p,q}$ with weights $0<p,q<1$ satisfying $p+q<1$ as in \eqref{eq:g_3}, and define $\beta_{p,q}:=(1-p)(1-q)(p+q)(p+q-p^2-q^2-pq)$. Then the probabilities $\PP(N_{p,q}=n)$ for $n\in\{3,4,5,6\}$ are given as in Table \ref{tab:thm}.
    \begin{table}[b!]
	\centering
	\renewcommand{\arraystretch}{2}
	\begin{tabular}{p{1.5cm}|p{14cm}}
		 & \hspace{0.15cm} $\PP(N_{p,q}=n)$ \\
		\hline
		$n=3$ 
            & \hspace{0.15cm} { \small $ \beta_{p,q}^{-1} \, \Big[ \; 2pq(1-p)(1-q)(p+q)(1-p-q) \; \Big] $ } \\[10pt]
        $n=4$ 
            & \hspace{0.15cm} { \small $ \beta_{p,q}^{-1} \, \Big[ \; 6p^2q^2(p+q)^2 + 2pq(12pq+1) - 22p^2q^2(p+q) $ } \\[5pt]
            & \hspace{1.5cm} { \small $ - p^2(5p^2q-12pq+2p+9q-p^2-1) - q^2(5pq^2-12pq+2q+9p-q^2-1)  \; \Big] $ } \\[10pt]
		$n=5$ 
            & \hspace{0.15cm} { \small $ \beta_{p,q}^{-1} \, \Big[ \; 6p^2q^2(p+q)(1-p-q) - 2pq(1-p-q)(p^2+q^2) $} \\[5pt]
            & \hspace{1.5cm} {\small $+ 2pq(p+q)(1-p-q) - 8p^2q^2(1-p-q) \; \Big] $ } \\[10pt]
		$n=6$ 
            & \hspace{0.15cm} { \small $ \beta_{p,q}^{-1} \, \Big[ \; 2 p^2 q^2 (1-p-q)^2 \; \Big] $ } 
	\end{tabular}
	\caption{The probabilities $\PP(N_{p,q}=n)$ for $n\in\{3,4,5,6\}$.}
	\label{tab:thm}
    \end{table}
\end{thm}

We will see in the course of this paper that the choice $p=q=1/3$ plays a special role. In fact, for these particular weights the probabilities $\PP(N_{p,q}=n)$ are maximized if $n\in\{3,5,6\}$ and minimized for $n=4$, see Lemmas \ref{lem:proportion_triangles}--\ref{lem:proportion_hexagons} below. The probabilities $\PP(N_{1/3,1/3}=n)$ are summarized in the following table: 
% \vspace{-0.25cm}
\begin{center}
		\renewcommand{\arraystretch}{2}
		\begin{tabular}{c|c|c|c|c c}
		& $n=3$ & $n=4$ & $n=5$ & $n=6$ &\\
		\cline{1-5}
		$\PP(N_{1/3,1/3}=n)$ & {\Large $2\over 9$} &  {\Large $7\over 12$} &  {\Large $1\over 6$} &  {\Large $1\over 36$} & .
	\end{tabular}
\end{center}
\vspace{0.14cm}
As a direct consequence of Theorem \ref{thm:main_theorem}, it is not difficult to confirm that $\EE N_{p,q}=4$. In addition, we can also derive the variance of the random variable $N_{p,q}$: 
 \begin{align*}
	\var N_{p,q} = \frac{4 p q (1-p-q)}{(1-p) (1-q) (p+q)}.
\end{align*}
We remark that $\var N_{p,q}$ takes its maximal value $1/2$ precisely if $p=q=1/3$.

% ====================================================================
\section{Preliminaries}
\label{sec:preliminaries}
% ====================================================================

Fix weights $0<p,q<1$ with $p+q<1$ and consider a Poisson line tessellation $X_{p,q}$ with directional distribution $G_{p,q}$ as in \eqref{eq:g_3}. Also recall from \eqref{eq:TypicalCell} the definition of the typical cell $\tc_{p,q}$ of $X_{p,q}$. A method to sample a random polygon with the same distribution as $\tc_{p,q}$ has been proposed in \cite{george_1987} and turns out to be rather powerful for our purpose.

\begin{figure}
	\centering
    \begin{subfigure}{0.5\textwidth}
    
        \centering
        \begin{tikzpicture}
            % GRID
    		% \draw[step=0.5cm,gray,very thin] (-3,-1) grid (3.5,5);

            % horizontal
            \draw [dashed] (-2,0)--(3.5,0) node[right, below, xshift=-2em]{\scriptsize horizontal};
    
            % Fünfeck zeichnen
            \coordinate (v1) at (0,0);
            \draw[-, name path=one] (v1)--++(120:3cm) node[midway, left, yshift=-4pt]{\small$z_1$} node[midway,style={sloped,allow upside down}]{\midarrow};
            \draw[-, name path=two] (v1)++(120:3cm)-- ++(60:2.5cm) node[midway, left, xshift=-2pt, yshift=2pt]{\small$z_2$} node[midway,style={sloped,allow upside down}]{\midarrow};
            \draw[-, name path=three] (v1)++(120:3cm)++(60:2.5cm)-- ++(0:1.5cm) node[midway, above, yshift=2pt]{\small$z_3$} node[midway,style={sloped,allow upside down}]{\midarrow};
            \draw[-, name path=four] (v1)++(120:3cm)++(60:2.5cm)++(0:1.5cm) -- ++(-60:1.5cm)node[midway, left, xshift=-2pt, yshift=-1pt]{\small$z_4$} node[midway,style={sloped,allow upside down}]{\midarrow};
            \draw[-, name path=five] (v1)++(120:3cm)++(60:2.5cm)++(0:1.5cm) ++(-60:1.5cm) -- (v1)node[midway, right, xshift=3pt, yshift=-2pt]{\small$z_5$} node[midway,style={sloped,allow upside down}]{\midarrow};

            % Schnittpunkte definieren
            \path [name intersections={of=one and five, by={A}}];
            \path [name intersections={of=one and two, by={B}}];
            \path [name intersections={of=two and three, by={C}}];
            \path [name intersections={of=three and four, by={D}}];
            \path [name intersections={of=four and five, by={E}}];
    
            % PLT andeuten, alle richtig außer E
            \draw[thin, gray] (A) -- ++(-60:1cm);
            \draw[thin, gray] (A) -- ++(-120:1cm);
            \draw[thin, gray] (B) -- ++(-120:1cm) node[gray, at end, xshift=-2pt, yshift=-4pt]{\small ${\ell_2}$};
            \draw[thin, gray] (B) -- ++(120:1cm) node[gray, at end, xshift=-4pt, yshift=6pt]{\small ${\ell_1}$};
            \draw[thin, gray] (C) -- ++(180:1cm) node[gray, at end, xshift=-6pt, yshift=0pt]{\small ${\ell_3}$};
            \draw[thin, gray] (C) -- ++(60:1cm);
            \draw[thin, gray] (D) -- ++(120:1cm);
            \draw[thin, gray] (E) -- ++(-60:1cm) node[gray, at end, xshift=5pt, yshift=-5pt]{\small ${\ell_4}$};
            \draw[thin, gray] (E) -- ++(60:1cm) node[gray, at end, xshift=5pt, yshift=5pt]{\small ${\ell_5}$};
    
            % Horizontalen an den Eckpunkten
            \draw[dashed] (B) -- ++(0:0.75cm);
            \draw[dashed] (C) -- ++(0:0.75cm);
            \draw[dashed] (D) -- ++(0:0.75cm);
            \draw[dashed] (E) -- ++(0:0.75cm);
    
            \filldraw[black] (A) circle circle (0.1pt) node[below, yshift=-3.4pt, xshift=1pt]{\small$v_1$};
            \filldraw[black] (B) circle circle (0.1pt) node[left]{\small$v_2$};
            \filldraw[black] (C) circle circle (0.1pt) node[left, yshift=6pt, xshift=2pt]{\small$v_3$};
            \filldraw[black] (D) circle circle (0.1pt) node[right, yshift=6pt, xshift=-2pt]{\small$v_4$};
            \filldraw[black] (E) circle circle (0.1pt) node[right, yshift=6pt, xshift=2pt]{\small$v_5$};

            % Winkel
            \draw (A) ++(0:.5) arc (0:60:.5) ; 
            \draw (A) ++(0:1.2) arc (0:120:1.2);
            \draw (B) ++(0:.5) arc (0:60:.5); 
            \draw (D) ++(0:0.5) arc (0:-60:.5); 
            \draw (E) ++(0:.5) arc (0:-120:.5);
    
            % Winkelbezeichnungen
            \draw node at (0.75,0.25){\small$\varphi_0$};
            \draw node at (1.5,0.25){\small$\varphi_1$};
            \draw node at (-0.7,2.85){\small$\varphi_2$};
            \draw node at (2.1,4.5){\small$\varphi_4$};
            \draw node at (2.8,3.2){\small$\varphi_5$};        
        \end{tikzpicture}

        \caption{Sequential construction of a polygon in cyclic order.}
        \label{fig:cyclic_constr}
    
    \end{subfigure}%
    \begin{subfigure}{0.5\textwidth}
        \centering
        \begin{tikzpicture}
            \draw [dashed, name path=ho] (-2.5,2)--(3,2) node[right, below, xshift=-4pt]{\scriptsize horizontal};

            \draw[-, thick, name path=one] (-2,1)--++(60:3.5cm) node[left, yshift=0pt]{$\ell_1$} node[midway,style={sloped,allow upside down}]{\midarrow};
            \draw[-, thick, name path=two] (1.5,3)--++(-120:3.5cm) node[left, yshift=0pt]{$\ell_2$} node[midway,style={sloped,allow upside down}]{\midarrow};

            % Schnittpunkte definieren
            \path [name intersections={of=one and ho, by={A}}];
            \path [name intersections={of=two and ho, by={B}}];
            
            % Namen Schnittpunkte
            \draw node at (A) [xshift=-8pt, yshift=8pt]{$v_1$};
            \draw node at (B) [xshift=-8pt, yshift=8pt]{$v_2$};
            
            % Winkel
            \draw (A) ++(0:1) arc (0:60:1) node[midway, xshift=7pt, yshift=7pt]{$\varphi_1$}; 
            \draw (B) ++(0:1) arc (0:-120:1) node[midway, xshift=7pt, yshift=-7pt]{$\varphi_2$};            
        \end{tikzpicture}
    \caption{Orientation of the angles}
    \label{fig:orientation_angles}
    \end{subfigure}%
    \caption{{Visualization of the concepts used for constructing polygons}
    }
\end{figure}
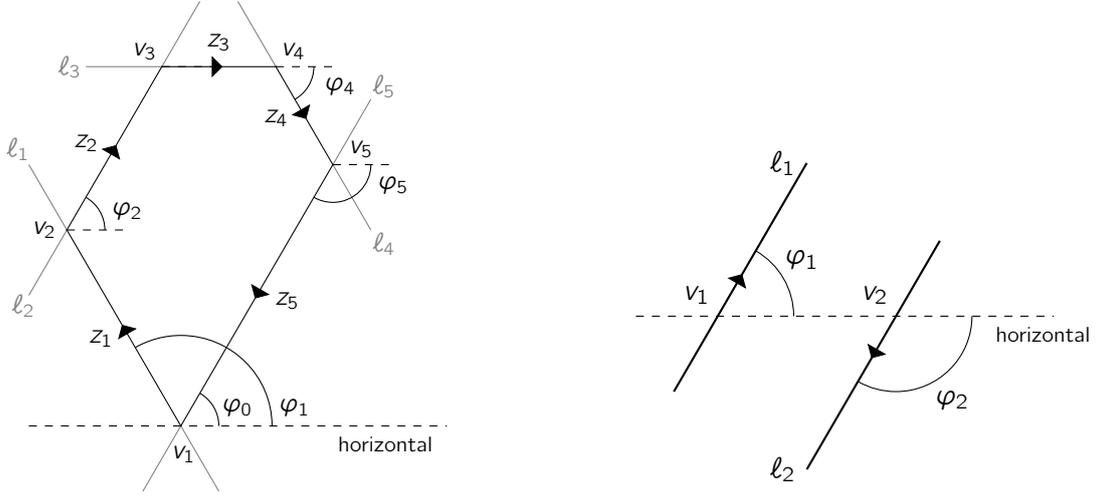

To explain it, fix $n\in\{3,4,5,6\}$ and note that an $n$-sided polygon is determined by the $n$ oriented lines $\ell_1,\ldots,\ell_n$ that support its sides which we think of being arranged in cyclic order as shown in Figure \ref{fig:cyclic_constr}. Alternatively, an $n$-sided polygon is determined by the following parameters: 
\begin{itemize}
	\item[(i)] the lengths $z_1,\ldots,z_{n}>0$ of the polygon's sides which are located on the lines $\ell_1,\ldots,\ell_{n}$,
	\item[(ii)] the angles $\varphi_0,\ldots,\varphi_{n-1}\in(-\pi,\pi)$, where $\varphi_i$ is the orientated angle at vertex $i$ that $\ell_i$ encloses with the eastern horizontal axis and where the sign of $\varphi_i$ is determined as explained in Figure~\ref{fig:orientation_angles}.
\end{itemize}
Putting $\varphi_n:=\varphi_0-\pi$, we next observe that the last two side lengths $z_{n-1}$ and $z_n$ can be determined from the remaining parameters because of the relation
\begin{equation}\label{eq:RemainingSides}
\sum_{i=1}^n z_i\sin\varphi_i = \sum_{i=1}^n z_i\cos\varphi_i=0,
\end{equation}
see \cite[Equation~(2.7)]{george_1987}. In what follows we shall write $\mathsf{poly}_n$ for the space of $n$-sided polygons in the plane whose lexicographically smallest vertex has coordinates $(0,0)$ and $P(\varphi_0,\ldots,\varphi_{n-1},z_1,\ldots,z_{n-2})\in\mathsf{poly}_n$ for the $n$-sided polygon determined by $\varphi_0,\ldots,\varphi_{n-1},z_1,\ldots,z_{n-2}$.

Applying this parametrization to the typical cell of the Poisson line tessellation $X_{p,q}$ and writing $Z_1,\ldots,Z_{n-2}$ for the random side lengths and $\Phi_0,\ldots,\Phi_{n-1}$ for the random orientation angles, a special case of the main result of \cite{george_1987} yields the following joint density for the random vector $(Z_1,\ldots,Z_{n-2},\Phi_0,\ldots\newline\ldots,\Phi_{n-1})$, see \cite[Equation (4.6)]{george_1987}. In this paper we use the convention that $G(\{\varphi\})=G(\{\pi-|\varphi|\})$ if $\varphi<0$.

\begin{lem}
\label{lem:density_z^n}
	Consider a Poisson line tessellation with directional distribution $G_{p,q}$ having weights $0<p,q<1$ which satisfy $p+q<1$, and fix $n\in\{3,4,5,6\}$. For  $\varphi\in\{0,\pi/3,2\pi/3\}$ define
    \begin{equation*}
        \lambda(\varphi) := p\vert\sin(\varphi)\vert + q\vert\sin(\pi/3-\varphi)\vert + (1-p-q)\vert\sin(2\pi/3-\varphi)\vert.
    \end{equation*}
    and $\lambda:=\sqrt{3}(p+q-p^2-q^2-pq)$. Then the conditional distribution of $\tc_{p,q}$ given $N_{p,q}=n$ is described by the random vector $(Z_1,\ldots,Z_{n-2},\Phi_0,\ldots,\Phi_{n-1})$ whose joint density with respect to the product of the Lebesgue measure on $(0,\infty)^{n-2}$ and $G_{p,q}^{\otimes n}$ is given by
    \begin{align}
	\label{eq:density_z^n}
		(\varphi_0,\ldots,\varphi_{n-1},z_1,\ldots,z_{n-2}) \, \longmapsto \; &
        \frac{2}{\lambda}\left(\frac{\sqrt{3}}{2}\right)^{n-1}
		\exp\Big( -\frac{1}{2}\sum_{i=1}^n z_i\lambda(\varphi_i) \Big) \nonumber \\[5pt]
        &\quad \times \one \big\{ P(\varphi_0,\ldots,\varphi_{n-1},z_1,\ldots,z_{n-2}) \in \mathsf{poly}_n \big\}.
	\end{align}
\end{lem}

One might think that integration of the density in  \eqref{eq:density_z^n} does not involve much effort. However, it turns out to be a rather subtle task, which at the same time requires some corrections of the method described in \cite{george_1987}. The difficulties arise from the indicator function of the event that the particular sequence of side lengths and orientation angles does indeed lead to an $n$-sided polygon. While integration with respect to the orientation angles is straight forward, attention has to be paid to the integration with respect to the side lengths. Correcting Equation~(3.12) in \cite{george_1987} gives for $i\in\{1,\ldots,n-2\}$ the upper limit of integration $\upperu_i$ for the variable $z_i$:
\begin{equation}
	\label{eq:u_i}
	\upperu_i:=\begin{cases}
		        -\csc(\varphi_i-\varphi_0) \sum\limits_{j=1}^{i-1} z_j\sin(\varphi_j-\varphi_0) &:\varphi_i<\varphi_0 \\
		        \infty &:\varphi_i\ge\varphi_0.
		    \end{cases}
\end{equation}
The lower limit of integration $\loweru_i$ for the variable $z_i$ is implicitly assumed to be zero in \cite{george_1987}. However, it will become clear from our computations that this is not always correct. Indeed, in particular situations some sides must have a minimum length strictly larger than zero due to given angles. Since these occurrences are highly dependent on the construction of the specific polygon, there seems to be no close form representation for these lower limits of integration. They will therefore be discussed in detail whenever they appear in what follows.

% ====================================================================
\section{Proof of Theorem \ref{thm:main_theorem}}
\label{sec:proofThm1}
% ====================================================================

% ====================================================================
\subsection{The triangle case}
\label{sec:triangles}
% ====================================================================

The probability that the typical cell $\tc_{p,q}$ is a triangle as been determined in  \cite[Theorem~1.1]{heerten_krecklenberg_thaele_2022}. To keep this paper self-contained, we briefly discuss the result and its derivation. Taking $n=3$ in the density in \eqref{eq:density_z^n} gives 
\begin{align*}
	(\varphi_0,\varphi_1,\varphi_2,z_1) \mapsto \;& \frac{3}{2\lambda}\exp\Big(-\frac12 \big( z_1\lambda(\varphi_1)+z_2\lambda(\varphi_2)+z_3\lambda(\varphi_3)\big)\Big)\one\big\{ P(\varphi_0,\varphi_1,\varphi_2,z_1)\in\mathsf{poly}_3\big\} .
\end{align*}
The difficult part in the integration of this density comes from the indicator function. For $n=3$ one has that
\begin{align*}
	\one\big\{ P(\varphi_0,\varphi_1,\varphi_2,z_1)\in\mathsf{poly}_3 \big\} &= \one\big\{ \varphi_0=0, \varphi_1=\pi/3, \varphi_2=2\pi/3,z_1\in(\loweru_{1},\upperu_{1}) \big\} \\
	&\qquad + \one\big\{\varphi_0=\pi/3, \varphi_1=2\pi/3, \varphi_2=0, z_1\in(\loweru_{2},\upperu_{2}) \big\} ,
\end{align*}
since it was argued in \cite{heerten_krecklenberg_thaele_2022} that there are only two different configurations $\triangle_1$, $\triangle_2$ of $(z_1,\varphi_0,\varphi_1,\varphi_2)$ that lead to a triangle. They are summarized in Table~\ref{tab:angles_for_triangles}. {In both cases, the lower integration} limit for $z_1$ is zero and the upper integration limit is given by $\infty$. Therefore, integrating the density of $(Z_1,\Phi_0,\Phi_1,\Phi_2)$ for $\triangle_1$ yields
\begin{align*}
	\int_{[0,\pi)^3}&\int_0^\infty \frac{3}{2\lambda} \, \exp\Big(-\frac12 \big( z_1\lambda(\varphi_1)+z_2\lambda(\varphi_2)+z_3\lambda(\varphi_3)\big)\Big) \\[5pt]
	&\hspace{3cm}\times \one\big\{\varphi_0=0, \varphi_1=\pi/3, \varphi_2=2\pi/3\big\} \,\dint z_1G^{\otimes 3}(\dint(\varphi_0, \varphi_1, \varphi_2)) \\[10pt]
	&= \frac{3pq(1-p-q)}{2\lambda}  \int_0^\infty \exp\Big(-\frac{z_1}{2} \big( \lambda(\varphi_1)+\lambda(\varphi_2)+\lambda(\varphi_3)\big)\Big)\, \dint z_1 \\[10pt]
	&= \frac{p q (1-p-q)}{p+q-p^2-q^2-pq}.
\end{align*}
Here, we used the fact that under the condition that $\varphi_0=0, \varphi_1=\pi/3, \varphi_2=2\pi/3$ the triangle has to be regular. In a similar way, the same result for $\triangle_2$ is obtained. Combining both cases, we recover  \cite[Theorem~1.1]{heerten_krecklenberg_thaele_2022} and have thus proved the first row of Table \ref{tab:thm}. Our findings are summarized in the following lemma, which also involves a discussion of extremal cases. The probability $\PP(N_{p,q}=3)$ is visualized in Figure~\ref{fig:plot_of_p3}.

\begin{table}[t!]
	\centering
	\renewcommand{\arraystretch}{1.5} 
	\begin{tabular}{c|c|c|c|c}
		Case & $\varphi_0$ & $\varphi_1$ & $z_1$ & $\varphi_2$  \\
		\hline
		$\triangle_1$ & $0$ & $\pi/3$ & $(0,\infty)$ & $2\pi/3$ \\
		$\triangle_2$ & $\pi/3$ & $2\pi/3$ & $(0,\infty)$ & $0$ 
	\end{tabular}
	\renewcommand{\arraystretch}{1} %back to normal 
	\caption{Values of angles $\varphi_i$ and intervals $(\loweru_1,\upperu_1)$ for side length $z_1$ resulting in a triangle.}
	\label{tab:angles_for_triangles}
\end{table}%\\

\begin{lem}
\label{lem:proportion_triangles}
    In the setup of Theorem~\ref{thm:main_theorem} it holds that 
    \begin{align*}
    % \label{eq:p_3}
        \PP(N_{p,q}=3)=\beta_{p,q}^{-1} \, \big[ \, 2pq(1-p)(1-q)(p+q)(1-p-q) \, \big].
    \end{align*}
    The maximum value for $\PP(N_{p,q}=3)$ is attained precisely if $p=q=1/3$ and is given by
    \begin{equation*}
        \max_{0<p+q<1}\PP(N_{p,q}=3)=\PP(N_{1/3,1/3}=3)=\frac29.
    \end{equation*}
\end{lem}

% ====================================================================
\subsection{The quadrilateral case}
\label{sec:quadrilaterals}
% ====================================================================

Since in comparison to the triangle case discussed above the results of this and the subsequent sections are new, we will discuss them in more detail. We start with the observation that the collection of quadrilaterals arising in the Poisson line tessellation $X_{p,q}$ can be subdivided into two classes: parallelograms ($\mathsf{para}$) and trapezoids ($\mathsf{trap}$). The method described in Section~\ref{sec:preliminaries} yields three possible angle configurations for the typical cell $\tc_{p,q}$ that belongs to $\mathsf{para}$ and six configurations leading to a quadrilateral in $\mathsf{trap}$, which are summarized in Table~\ref{tab:angles_for_quadrilaterals}. 

\begin{figure}[t!]
    \centering
    \begin{tikzpicture}
		\draw [dashed] (-1,0)--(4.5,0) node[right, below, xshift=-2em]{\scriptsize horizontal};

        \coordinate (v1) at (0,0);
        \draw[-, name path=one] (v1)--++(120:2.5cm) node[midway, left]{$z_1$};
        \draw[-, name path=twoa] (v1)++(120:2.5cm)-- ++(0:3.5cm) node[midway, above]{$z_2$};
        \draw[-, name path=twob] (v1)++(120:2.5cm)++(0:3.5cm)-- ++(0:2.5cm);
        \draw[-, name path=three] (v1)++(120:2.5cm)++(0:6cm)-- ++(-120:2.5cm);
        \draw[-, name path=four] (v1)++(120:2.5cm)++(0:6cm)-- ++(-120:2.5cm) -- (v1);
        \path[name path=four] (v1)++(120:2.5cm)++(0:6cm)++(-120:2.5cm)--(v1);
        
        \draw[dashed, line width=0.4mm, lightgray, name path=d] (v1)++(-120:0.5cm)--++(60:3.5cm);
        \draw[line width=0.4mm, lightgray, name path=d2] (1,0)++(-120:0.5cm)--++(60:3.5cm) node[right, xshift=-2pt]{$\ell_3$};

        \path [name intersections={of=one and twoa, by={B}}];
        \path [name intersections={of=twob and three, by={C}}];
        \path [name intersections={of=three and four, by={D}}];
        \path [name intersections={of=four and one, by={A}}];
        \path [name intersections={of=twob and d2, by={W}}];

        \filldraw[black] (A) circle circle (0.1pt) node[below, xshift=3pt]{$v_1$};
        \filldraw[black] (B) circle circle (0.1pt) node[left]{$v_2$};
        \filldraw[black] (W) circle circle (0.1pt) node[right, xshift=-2pt, yshift=-8pt]{$v_3$};
    \end{tikzpicture}
    \caption{Minimal length of $z_2$ in $\squad_{\,4}$.}
    \label{fig:quadrilaterals_case4}
\end{figure}
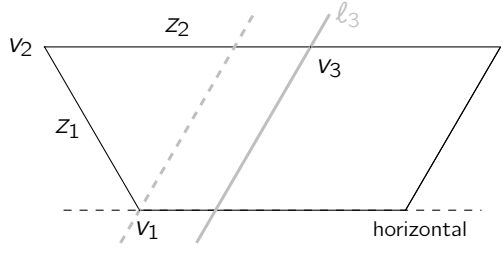

Taking $n=4$ in Lemma~\ref{lem:density_z^n} implies that the joint density of $(\Phi_0,\Phi_1,\Phi_2,\Phi_3,Z_1,Z_2)$ is given by
\begin{equation}
\label{eq:density_z4}
    (\varphi_0,\dots,\varphi_3,z_1,z_2) \longmapsto
    \frac{3\sqrt{3}}{4\lambda}\, \exp\left(-\frac{1}{2} \sum_{i=1}^4 z_i\lambda(\varphi_i)\right)\one\big\{ P(\varphi_0,\varphi_1,\varphi_2,\varphi_3,z_1,z_2)\in\mathsf{poly}_4\big\} .
\end{equation}
We express the remaining two side lengths using \eqref{eq:RemainingSides}. This leads to
\begin{equation*}
    z_3=
    \begin{cases}
        z_1&:\text{in cases } 1, 2, 4, 5, 7, 8 \\
        z_1-z_2&:\text{in cases } 3, 9 \\
        z_1+z_2&:\text{in case\phantom{s} } 6,
    \end{cases}
    \qquad
    z_4=
    \begin{cases}
        z_2&:\text{in cases } 1, 3, 5, 6, 8, 9\\
        z_1+z_2&:\text{in cases } 2, 7 \\
        -z_1+z_2&:\text{in case\phantom{s} } 4.
    \end{cases}
\end{equation*}
In order to integrate the density in \eqref{eq:density_z4} with respect to $z_1$ and $z_2$, the upper integration limits $\overline{u}_1$ and $\overline{u}_2$ in \eqref{eq:u_i} need to be determined. This yields $\overline{u}_1=\infty$, independently of the individual case, and
\begin{equation*}
    % u_1=\infty \quad \text{for all cases and}
    % \quad 
    \overline{u}_2=
    \begin{cases}
        z_1&:\text{in cases } 3,9 \\
        \infty&:\text{else.}%\text{in cases } 1,2,4,5,6,7,8. \\
    \end{cases}
\end{equation*}
As already discussed in Section~\ref{sec:preliminaries}, there is a need to consider lower integration limits different from zero. Indeed, for quadrilaterals this situation appears precisely for the trapezoid described by line $4$ in Table \ref{tab:angles_for_quadrilaterals}, from here on denoted by $\squad_{\,4}$. As this acts as a model case for later purposes (when $n=5,6$), we discuss this issue here in detail. {
Therefore, we assume that $(\varphi_0,\varphi_1,\varphi_2,\varphi_3)=(0,2\pi/3,0,-2\pi/3)$. Depending on the value of $z_2$, it is possible that the line $\ell_3$, on which the polygon side with length $z_3$ is located, intersects the first polygon side with length $z_1$ prior to intersecting the horizontal line. This is illustrated in Figure~\ref{fig:quadrilaterals_case4}, which shows that the line $\ell_3$ cannot be on the left side of the parallel dashed line, since otherwise the construction leads to a triangle. Therefore, $z_2$ must be at least $z_1$, since the two sides with length $z_1$, $z_2$ together with the dashed line comprise a regular triangle. In other words $\underline{u}_2=z_1$.
}

\begin{table}[t!]
	\centering
	\renewcommand{\arraystretch}{1.5} 
	\begin{tabular}{c|c|c|c|c|c|c|c}
		Case & $\varphi_0$ & $\varphi_1$ & $z_1$ %(\underline{u}_1,\overline{u}_1)$ 
        & $\varphi_2$ & $z_2$ & $\varphi_3$ & type of quadrilateral\\
		\hline
		$\squad_{\,1}$ & $0$ & $\pi/3$ & $(0,\infty)$ & $0$ & $(0,\infty)$ & $-2\pi/3$ & $\mathsf{para}$ \\
		$\squad_{\,2}$ & $0$ & $\pi/3$ & $(0,\infty)$ & $0$ & $(0,\infty)$ & $-\pi/3$ & $\mathsf{trap}$ \\
		$\squad_{\,3}$ & $0$ & $\pi/3$ & $(0,\infty)$ & $-\pi/3$ & $(0,z_1)$ & $-2\pi/3$ & $\mathsf{trap}$\\
		$\squad_{\,4}$ & $0$ & $2\pi/3$ & $(0,\infty)$ & $0$ & $(z_1,\infty)$ & $-2\pi/3$ & $\mathsf{trap}$ \\
		$\squad_{\,5}$ & $0$ & $2\pi/3$ & $(0,\infty)$ & $0$ & $(0,\infty)$ & $-\pi/3$ & $\mathsf{para}$ \\
		$\squad_{\,6}$ & $0$ & $2\pi/3$ & $(0,\infty)$ & $\pi/3$ & $(0,\infty)$ & $-\pi/3$ & $\mathsf{trap}$ \\
		$\squad_{\,7}$ & $\pi/3$ & $2\pi/3$ & $(0,\infty)$ & $\pi/3$ & $(0,\infty)$ & $0$ & $\mathsf{trap}$ \\
		$\squad_{\,8}$ & $\pi/3$ & $2\pi/3$ & $(0,\infty)$ & $\pi/3$ & $(0,\infty)$ & $-\pi/3$ & $\mathsf{para}$ \\
		$\squad_{\,9}$ & $\pi/3$ & $2\pi/3$ & $(0,\infty)$ & $0$ & $(0,z_1)$ & $-\pi/3$ & $\mathsf{trap}$ \\
	\end{tabular}
	\renewcommand{\arraystretch}{1} %back to normal 
	\caption{Values of angles $\varphi_i$ and intervals $(\loweru_i,\upperu_i)$ for side lengths $z_i$ resulting in a parallelogram ($\mathsf{para}$) or a trapezoid ($\mathsf{trap}$).}
	\label{tab:angles_for_quadrilaterals}
\end{table}%

\begin{figure}[b!]
	\begin{subfigure}{0.5\textwidth}
		\centering
		\includegraphics[scale=0.6]{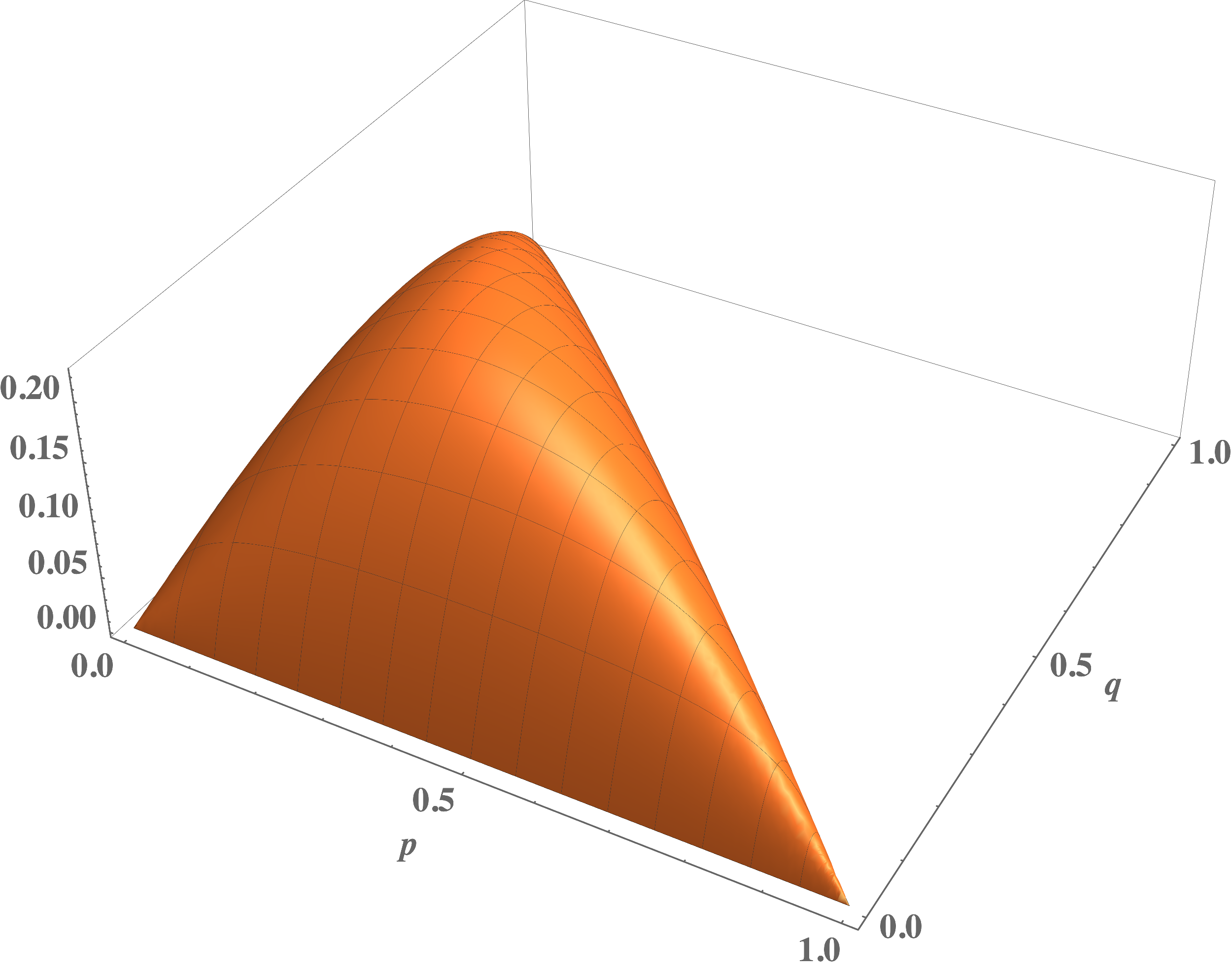}
		\caption{Plot of $\PP(N_{p,q}=3)$ as in Lemma~\ref{lem:proportion_triangles}.}
		\label{fig:plot_of_p3}
	\end{subfigure}%
	\begin{subfigure}{0.5\textwidth}
		\centering
		\includegraphics[scale=0.6]{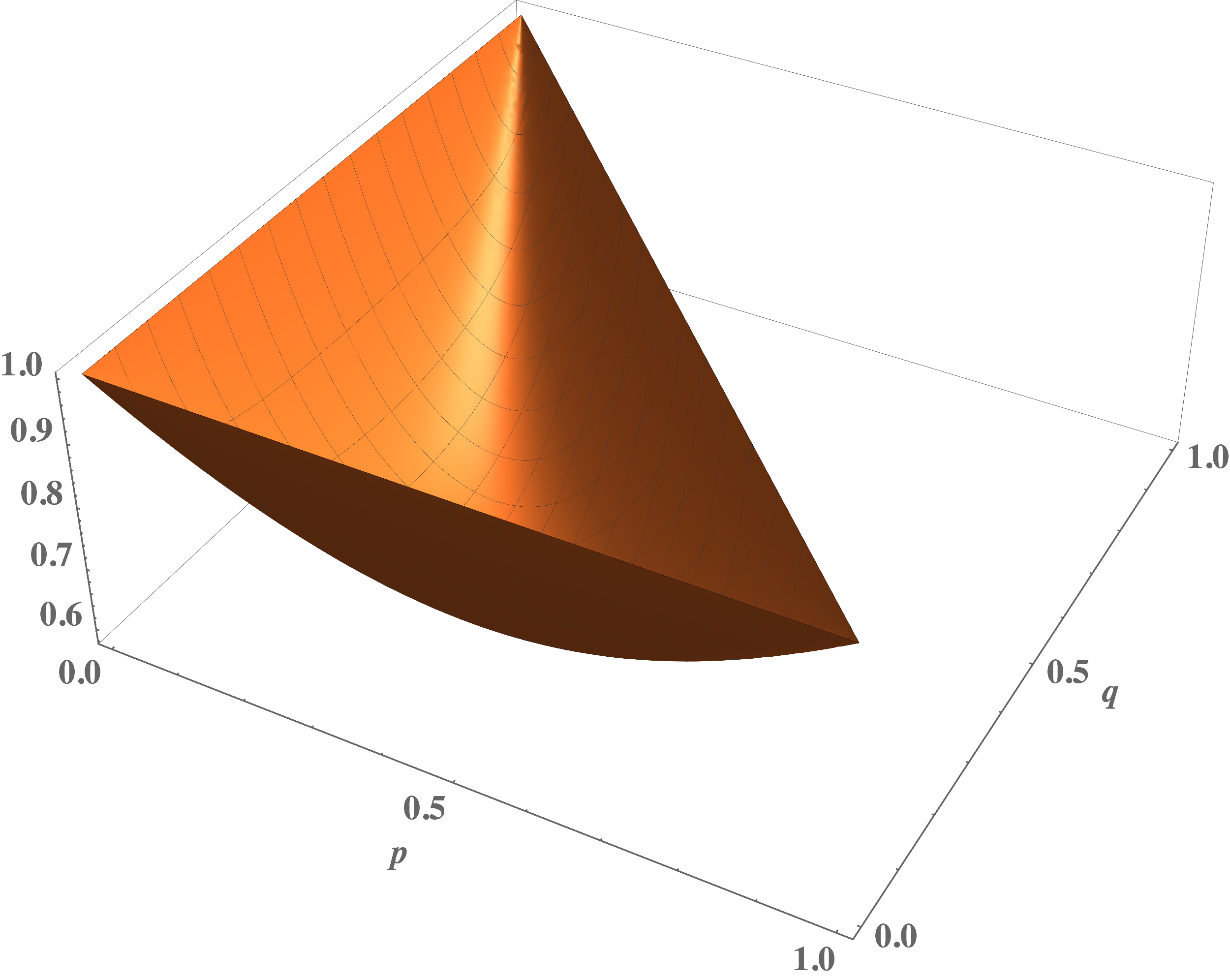}
		\caption{Plot of $\PP(N_{p,q}=4)$ as in Lemma~\ref{lem:proportion_quadrilaterals}.}
		\label{fig:plot_of_p4}
	\end{subfigure}%
	\caption{Plots of $\PP(N_{p,q}=3)$ and $\PP(N_{p,q}=4)$.}
\end{figure}

Taking into account the above considerations, integration of the density in \eqref{eq:density_z4} yields
\begin{align*}
    \int_{[0,\pi)^4}&\int_{0}^{\infty} \int_{z_1}^{\infty} \frac{3\sqrt{3}}{4\lambda}\exp\Big( -\frac{\sqrt{3}}{2} \big(p (z_1-z_2)+z_2\big)\Big)  \\[5pt]
    &\hspace{1.5cm}\times \one\big\{\varphi_0=0, \varphi_1=2\pi/3, \varphi_2=0,\varphi_3=-2\pi/3 \big\} \,\dint z_2\dint z_1G^{\otimes 4}(\dint(\varphi_0,\varphi_1,\varphi_2,\varphi_3)) \\[10pt]
    &= \frac{3\sqrt{3}p^2q(1-p-q)}{4\lambda} \int_{0}^{\infty} \int_{z_1}^{\infty} \exp\Big( -\frac{\sqrt{3}}{2} \big(p (z_1-z_2)+z_2\big)\Big) \dint z_2\dint z_1 \\[10pt]
    &= \frac{3p^2q(1-p-q)}{2\lambda(p-1)} \int_{0}^{\infty} -e^{-\frac{\sqrt{3}}{2}\,z_1}\, \dint z_1 \\[10pt]
    &= \frac{\sqrt{3}p^2q(1-p-q)}{\lambda(p-1)},
\end{align*}
where in the first step we used that $G(\{0\})=p$, $G(\{2\pi/3\})=1-p-q$ and $G(\{-2\pi/3\})=G(\{\pi/3\})=q$, according to our convention of signs. Now, inserting the value for $\lambda$ from Lemma \ref{lem:density_z^n} yields
\begin{align*}
    \PP(\tc_{p,q}\in\squad_{\,4}) = \frac{p^2q(1-p-q)}{(1-p)\left(p+q-p^2-q^2-pq\right)}.
\end{align*}
The other eight cases can be dealt with in the same way, which leads to the following result, illustrated in Figure \ref{fig:plot_of_p4}.
\begin{lem}
\label{lem:proportion_quadrilaterals}
    In the setup of Theorem~\ref{thm:main_theorem} it holds that 
    \begin{align*}
        \PP(N_{p,q}=4) &= 
        	\beta_{p,q}^{-1}\, \Big[ \; 6p^2q^2(p+q)^2 + 2pq(12pq+1) - 22p^2q^2(p+q)   \\[5pt]
        	& - p^2(5p^2q-12pq+2p+9q-p^2-1) - q^2(5pq^2-12pq+2q+9p-q^2-1)  \; \Big].
    \end{align*}
    The minimal value for $\PP(N_{p,q}=4)$ is attained precisely if $p=q=1/3$ and is given by
\begin{equation*}
	\min_{0<p+q<1}\PP(N_{p,q}=4)=\PP(N_{1/3,1/3}=4)=\frac{7}{12}.
\end{equation*}
\end{lem}

\begin{figure}[t!]
    \centering
    \begin{subfigure}{0.5\textwidth}
        \centering
        \includegraphics[scale=0.5]{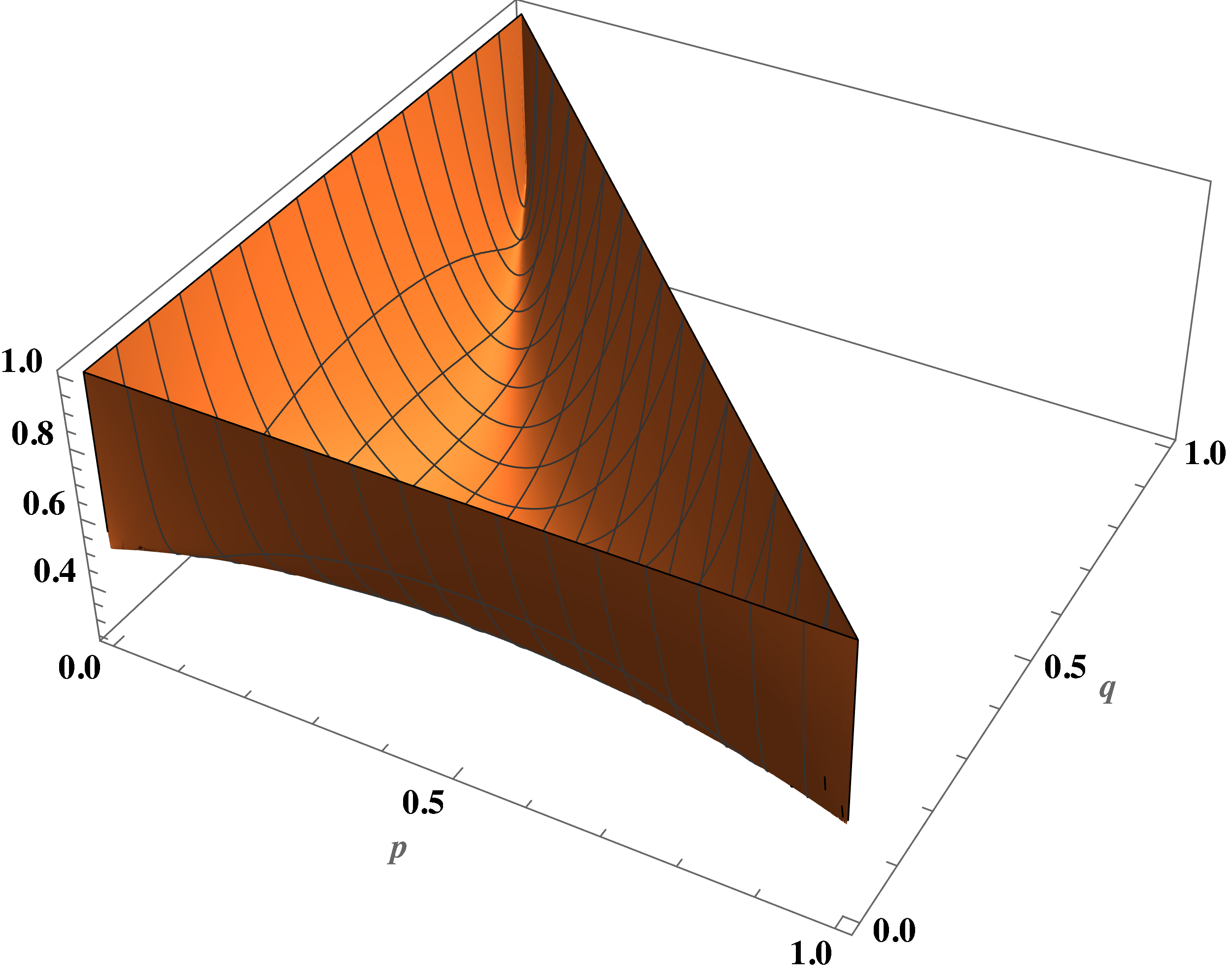}
        \caption{Plot of $\PP(\tc_{p,q}\in\mathsf{para})$.}
    \end{subfigure}%
    \begin{subfigure}{0.5\textwidth}
        \centering
        \includegraphics[scale=0.5]{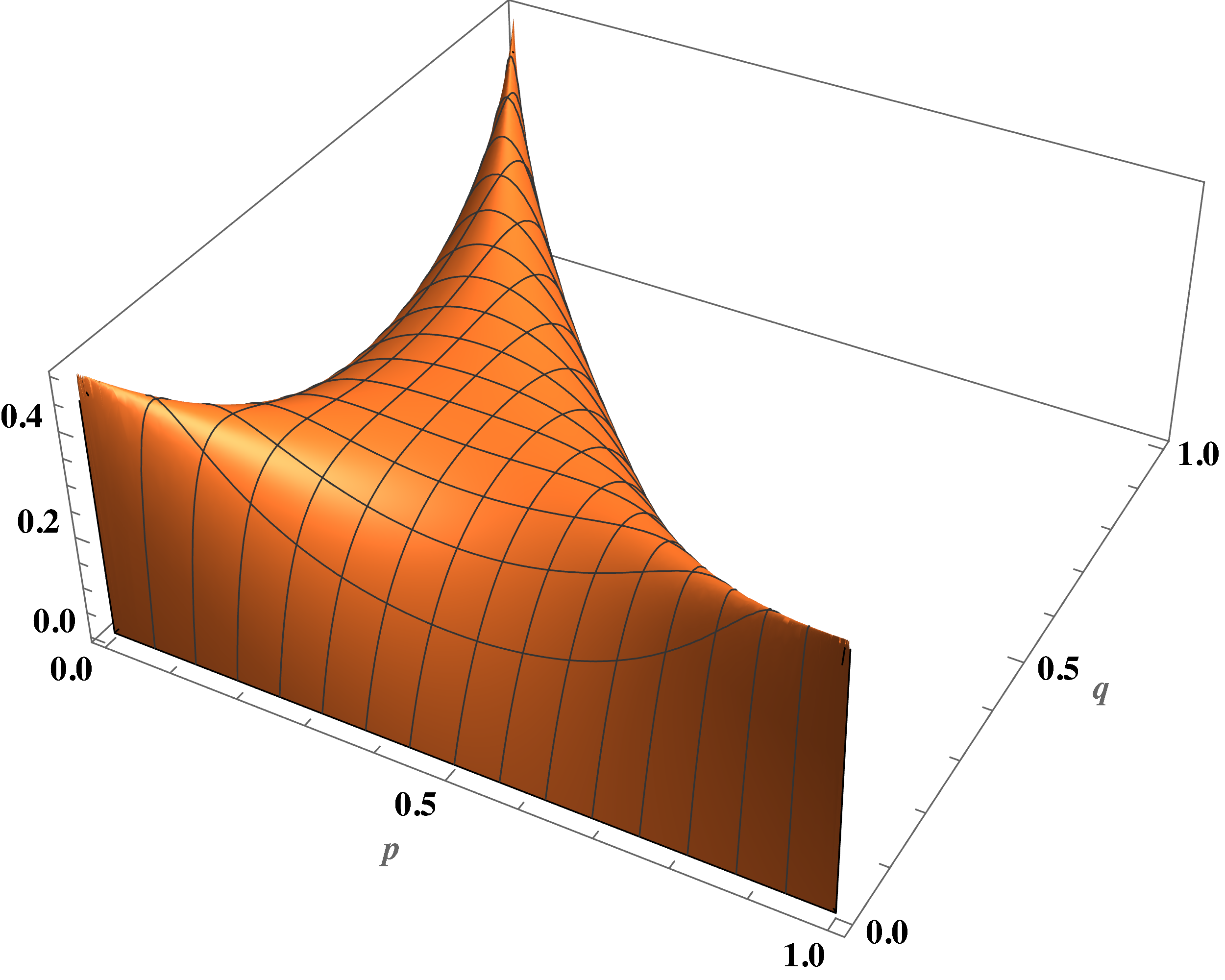}
        \caption{Plot of $\PP(\tc_{p,q}\in\mathsf{trap})$.}
    \end{subfigure}%
    \caption{Plots of $P(N_{p,q}=4)$ in the subcases of parallelograms and trapezoids.}
    \label{fig:plots_of_p4_subcases}
\end{figure}

\begin{remark}
\label{rem:proportion_trapezoids_parallelograms}
    It is be possible to refine Lemma \ref{lem:proportion_quadrilaterals} in the following way. Writing $\PP(N_{p,q}=4)=\PP(\tc_{p,q}\in\mathsf{para})+\PP(\tc_{p,q}\in\mathsf{trap})$, we have that
    \begin{align*}
        \PP(\tc_{p,q}\in\mathsf{para}) &= \beta_{p,q}^{-1} \, \Big[ \, p^4(1-q)-2p^3(q-1)^2+q^2(1-p)(q-1)^2 \\[5pt] 
         &\hspace{5cm} +2pq^2(p-1)+p^2(-2 q^3+6 q^2-3 q+1) \,\Big], \\[10pt]
        \PP(\tc_{p,q}\in\mathsf{trap}) &= \beta_{p,q}^{-1} \, \Big[ \, -5pq(1-p-q) -2(p+q)(1-p)(1-q) -2pq \,\Big].
    \end{align*}
    	Both probabilities are visualized in \autoref{fig:plots_of_p4_subcases}. 
    Similar to the total probability $\PP(N_{p,q}=4)$, the parallelogram probability $\PP(\tc_{p,q}\in\mathsf{para})$ also has a minimum that is attained precisely if $p=q=1/3$ and is given by $1/4$. On the other hand, we observe a local maximum of the trapezoid probability $\PP(\tc_{p,q}\in\mathsf{trap})$ around $p=q=1/3$ of $1/3$.
\end{remark}

% ====================================================================
\subsection{The pentagon case}
\label{sec:pentagons}
% ====================================================================

In this section, we consider the probability that the typical cell $\tc_{p,q}$ has five vertices. We have to distinguish between six different types of configurations resulting in pentagons, see Table~\ref{tab:angles_for_pentagons}. Note that one case is subdivided into two subcases $\pentago_{\,2.1}$ and $\pentago_{\,2.2}$. This circumstance will be discussed later.
\begin{table}[b!]
    \centering
    \renewcommand{\arraystretch}{1.5} 
    \begin{tabular}{p{0.65cm}|c|c|c|c|c|c|c|c}
    & $\varphi_0$ & $\varphi_1$ & $z_1$ & $\varphi_2$ & $z_2$ & $\varphi_3$ & $z_3$ & $\varphi_4$ \\
    \hline
    $\pentago_{\,1}$ & $0$ & $\pi/3$ & $(0,\infty)$ & $0$ & $(0,\infty)$ & $-\pi/3$ & $(0,z_1)$ & $-2\pi/3$ \\[10pt]
        \begin{tabular}{@{}c@{}}
            $\pentago_{\,2.1}$ \\
            $\pentago_{\,2.2}$ \\
        \end{tabular} 
        & $0$ & $2\pi/3$ & 
        \begin{tabular}{@{}c@{}}
            $(0,\infty)$ \\
            $(z_2,\infty)$ \\
        \end{tabular}
        & $0$ & 
        \begin{tabular}{@{}c@{}}
            $(z_1,\infty)$ \\
            $(0,\infty)$ \\
        \end{tabular}
        & $-\pi/3$ & 
        \begin{tabular}{@{}c@{}}
            $(0,z_1)$ \\
            $(z_1-z_2,z_1)$ \\
        \end{tabular}
        & $-2\pi/3$ \\[20pt]
     $\pentago_{\,3}$ & $0$ & $2\pi/3$ & $(0,\infty)$ &$\pi/3$ & $(0,\infty)$ & $0$ & $(0,\infty)$ & $-\pi/3$ \\[10pt]
     $\pentago_{\,4}$ & $0$ & $2\pi/3$ & $(0,\infty)$ &$\pi/3$ & $(0,\infty)$ & $0$ & $(z_1,\infty)$ & $-2\pi/3$ \\[10pt]
     $\pentago_{\,5}$ & $0$ & $2\pi/3$ & $(0,\infty)$ & $\pi/3$ & $(0,\infty)$ & $-\pi/3$ & $(z_1,z_1+z_2)$ & $-2\pi/3$ \\[10pt]
     $\pentago_{\,6}$ & $\pi/3$ & $2\pi/3$ & $(0,\infty)$ & $\pi/3$ & $(0,\infty)$ & $0$ & $(0,z_1)$ & $-\pi/3$ \\[10pt]
    \end{tabular}
    \renewcommand{\arraystretch}{1} %back to normal 
    \caption{Values of angles $\varphi_i$ and intervals $(\loweru_i,\upperu_i)$ for side lengths $z_i$ resulting in pentagons.}
    \label{tab:angles_for_pentagons}
 \end{table}
Inserting $n=5$ into \eqref{eq:density_z^n} yields the joint density
\begin{align}
\label{eq:density_z5}
    (\varphi_0,\dots,\varphi_4,z_1,z_2,z_3) &\longmapsto
    \frac{9}{8\lambda}\, \exp\bigg(-\frac{1}{2} \sum_{i=1}^5 z_i\lambda(\varphi_i)\bigg) \nonumber \\[5pt]
    &\hspace{2cm} \times \one\big\{ P(\varphi_0,\varphi_1,\varphi_2,\varphi_3,\varphi_4,z_1,z_2,z_3)\in\mathsf{poly}_5\big\}
\end{align}
for the random vector $(\Phi_0,\ldots,\Phi_4,Z_1,Z_2,Z_3)$.
Similar to Section~\ref{sec:quadrilaterals}, we use \eqref{eq:RemainingSides} to express the remaining side lengths by means of the others. This leads to
\begin{equation*}
    z_4=
    \begin{cases}
        z_1-z_3,&\quad\text{in cases } 1,2,6  \\
        z_1+z_2,&\quad \text{in cases } 3,4 \\
        z_1+z_2-z_3,&\quad \text{in case\phantom{s} } 5,
    \end{cases}
    \qquad
    z_5=
    \begin{cases}
        z_2+z_3,&\quad\text{in cases } 1,3,6  \\
        -z_1+z_2+z_3,&\quad \text{in cases } 2 \\
        -z_1+z_3,&\quad \text{in case\phantom{s} } 4,5.
    \end{cases}
\end{equation*}
The upper integral limits $\upperu_i$ for $i=1,2,3$ can again be obtained by evoking \eqref{eq:u_i}:  $\upperu_1,\upperu_2=\infty$ in all cases and
\begin{equation}
\label{eq:u3_pentagons}
   \upperu_3=
    \begin{cases}
        z_1,&\quad\text{in cases } 1,2,6 \\
        \infty,&\quad\text{in case\phantom{s} } 3,4 \\
        z_1+z_2,&\quad\text{in cases } 5. 
    \end{cases}
\end{equation}
As explained in Section~\ref{sec:preliminaries}, there are cases that need special attention when it comes to the lower integral limits for some side lengths. Here, these are cases $2,4$ and $5$. For better readability, we adopt our prior notation and denote the six different cases of pentagons by $\pentago_{\,i}$ for $i=1,\dots,6$. In some of these cases, we observe similar issues as in the quadrilateral case discussed in the previous section, where the line coinciding with the penultimate side of the polygon (here, these are $\ell_4$ and $z_4$, respectively) must not intersect the first side with length $z_1$ prior to intersecting the horizontal. 

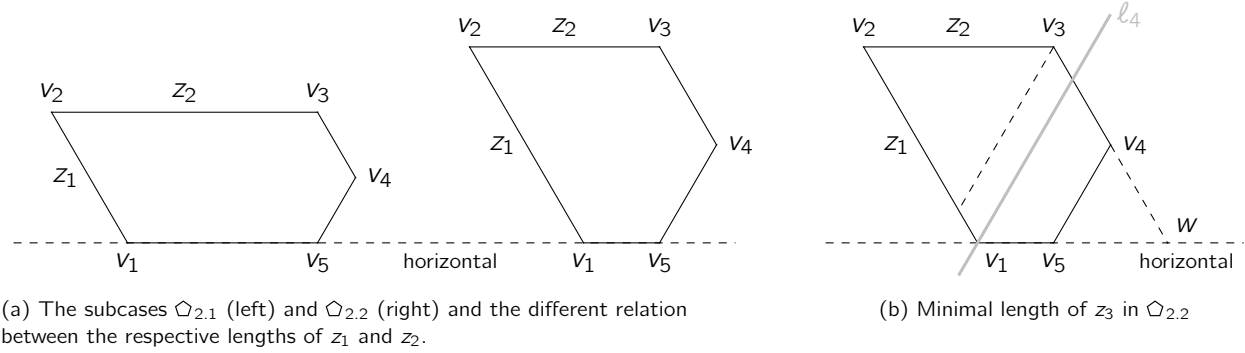
\begin{figure}[t!]
    \centering
    \begin{subfigure}{0.6\textwidth}
        \centering 
        \begin{tikzpicture}
    		\draw [dashed] (-2,0)--(7.5,0) node[midway, below, xshift=1cm]{\scriptsize horizontal};
    
            %case 2.1
            \coordinate (v1) at (-0.5,0);
            \draw[-, name path=one1] (v1)--++(120:2cm) node[midway, left]{$z_1$};
            \draw[-, name path=two1] (v1)++(120:2cm)-- ++(0:3.5cm) node[midway, above]{$z_2$};
            \draw[-, name path=three1] (v1)++(120:2cm)++(0:3.5cm)-- ++(-60:1cm);
            \draw[-, name path=four1] (v1)++(120:2cm)++(0:3.5cm)++(-60:1cm)-- ++(-120:1cm);
            \draw[-, name path=five1] (v1)++(120:2cm)++(0:3.5cm)++(-60:1cm)++(-120:1cm)-- (v1);
    
            \path [name intersections={of=one1 and two1, by={B1}}];
            \path [name intersections={of=two1 and three1, by={C1}}];
            \path [name intersections={of=three1 and four1, by={D1}}];
            \path [name intersections={of=four1 and five1, by={E1}}];
            \path [name intersections={of=five1 and one1, by={A1}}];
    
            \filldraw[black] (A1) circle (0.1pt) node[below]{$v_1$};
            \filldraw[black] (B1) circle (0.1pt) node[above]{$v_2$};
            \filldraw[black] (C1) circle (0.1pt) node[above]{$v_3$};
            \filldraw[black] (D1) circle (0.1pt) node[right]{$v_4$};
            \filldraw[black] (E1) circle (0.1pt) node[below]{$v_5$};
            
            %case 2.2
            \coordinate (vv1) at (5.5,0);
            \draw[-, name path=one2] (vv1)--++(120:3cm) node[midway, left]{$z_1$};
            \draw[-, name path=two2] (vv1)++(120:3cm)-- ++(0:2.5cm) node[midway, above]{$z_2$};
            \draw[-, name path=three2] (vv1)++(120:3cm)++(0:2.5cm)-- ++(-60:1.5cm);
            \draw[-, name path=four2] (vv1)++(120:3cm)++(0:2.5cm)++(-60:1.5cm)-- ++(-120:1.5cm);
            \draw[-, name path=five2] (vv1)++(120:3cm)++(0:2.5cm)++(-60:1.5cm)++(-120:1.5cm)-- (vv1);
            
            \path [name intersections={of=one2 and two2, by={B2}}];
            \path [name intersections={of=two2 and three2, by={C2}}];
            \path [name intersections={of=three2 and four2, by={D2}}];
            \path [name intersections={of=four2 and five2, by={E2}}];
            \path [name intersections={of=five2 and one2, by={A2}}];
            
            \filldraw[black] (A2) circle (0.1pt) node[below]{$v_1$};
            \filldraw[black] (B2) circle (0.1pt) node[above]{$v_2$};
            \filldraw[black] (C2) circle (0.1pt) node[above]{$v_3$};
            \filldraw[black] (D2) circle (0.1pt) node[right]{$v_4$};
            \filldraw[black] (E2) circle (0.1pt) node[below]{$v_5$};
        
        \end{tikzpicture}
    \caption{
        The subcases $\pentago_{\,2.1}$ (left) and $\pentago_{\,2.2}$ (right) and the different relation \\ between the respective lengths of $z_1$ and $z_2$.
    }
    \label{fig:pentagon_cases_2.1_2.2}
    \end{subfigure}%
    \begin{subfigure}{0.4\textwidth}
        \centering
        \begin{tikzpicture}
    		\draw [dashed, name path=horizontal] (-2,0)--(3.5,0) node[at end, below, xshift=-0.75cm]{\scriptsize horizontal};
            
            \coordinate (v1) at (0,0);
            \draw[-, name path=one2] (v1)--++(120:3cm) node[midway, left]{$z_1$};
            \draw[-, name path=two2] (v1)++(120:3cm)-- ++(0:2.5cm) node[midway, above]{$z_2$};
            \draw[-, name path=three2] (v1)++(120:3cm)++(0:2.5cm)-- ++(-60:1.5cm);
            \draw[-, name path=four2] (v1)++(120:3cm)++(0:2.5cm)++(-60:1.5cm)-- ++(-120:1.5cm);
            \draw[-, name path=five2] (v1)++(120:3cm)++(0:2.5cm)++(-60:1.5cm)++(-120:1.5cm)-- (v1);
            
            \path [name intersections={of=one2 and two2, by={B2}}];
            \path [name intersections={of=two2 and three2, by={C2}}];
            \path [name intersections={of=three2 and four2, by={D2}}];
            \path [name intersections={of=four2 and five2, by={E2}}];
            \path [name intersections={of=five2 and one2, by={A2}}];

            \filldraw[black] (A2) circle (0.1pt) node[below, xshift=0.25cm]{$v_1$};
            \filldraw[black] (B2) circle (0.1pt) node[above]{$v_2$};
            \filldraw[black] (C2) circle (0.1pt) node[above]{$v_3$};
            \filldraw[black] (D2) circle (0.1pt) node[right]{$v_4$};
            \filldraw[black] (E2) circle (0.1pt) node[below]{$v_5$};

            % dashed lines and supporting points 
            
            \draw[dashed, name path=d1] (D2)-- ++(-60:1.5cm);
            \path [name intersections={of=horizontal and d1, by={W1}}];
            \filldraw[black] (W1) circle (0.1pt) node[above, xshift=0.25cm]{$w$};

            \draw[dashed, name path=d2] (C2)-- ++(-120:2.5cm);
            % \draw[dashed, name path=d2] (A2)-- ++(60:2.5cm);

            \draw[-, line width=0.4mm, lightgray, name path=l4] (A2)++(-120:0.5cm)-- ++(60:4cm) node[right, xshift=-2pt]{$\ell_4$};
         \end{tikzpicture}
    \caption{Minimal length of $z_3$ in $\pentago_{\,2.2}$ \\ \phantom{x} }
    \label{fig:pent2.2}
    \end{subfigure}%
    
    \caption{
        The pentagon special case $\pentago_{\,2}$ with its two subcases $\pentago_{\,2.1}$ and $\pentago_{\,2.2}$.
    }
    \label{fig:PentCases2}
\end{figure}

We start with $\pentago_{\,2}$, which itself splits into two subcases, denoted by $\pentago_{\,2.1}$ and $\pentago_{\,2.2}$. This is due to the relation of $z_1$ and $z_2$, resulting in two similar pentagons, one more vertically stretched and the other more horizontally, see \autoref{fig:pentagon_cases_2.1_2.2} for an illustration.
\begin{figure}[b!]
    \centering
     \begin{subfigure}{0.6\textwidth}
        \centering
        \begin{tikzpicture}
    		\draw [dashed, name path=horizontal] (-1,0)--(4.5,0) node[at end, below, xshift=-0.75cm]{\scriptsize horizontal};
            
            \coordinate (v1) at (0,0);
            \draw[-, name path=one2] (v1)--++(120:1cm) node[midway, left]{$z_1$};
            \draw[-, name path=two2] (v1)++(120:1cm)-- ++(60:2cm) node[midway, left]{$z_2$};
            \draw[-, name path=three2] (v1)++(120:1cm)++(60:2cm)-- ++(0:3cm);
            \draw[-, name path=four2] (v1)++(120:1cm)++(60:2cm)++(0:3cm)-- ++(-120:3cm);
            \draw[-, name path=five2] (v1)++(120:1cm)++(60:2cm)++(0:3cm)++(-120:3cm)-- (v1);
            
            \path [name intersections={of=one2 and two2, by={B2}}];
            \path [name intersections={of=two2 and three2, by={C2}}];
            \path [name intersections={of=three2 and four2, by={D2}}];
            \path [name intersections={of=four2 and five2, by={E2}}];
            \path [name intersections={of=five2 and one2, by={A2}}];

            \filldraw[black] (A2) circle (0.1pt) node[below, xshift=0.25cm]{$v_1$};
            \filldraw[black] (B2) circle (0.1pt) node[left]{$v_2$};
            \filldraw[black] (C2) circle (0.1pt) node[above]{$v_3$};
            \filldraw[black] (D2) circle (0.1pt) node[above]{$v_4$};
            \filldraw[black] (E2) circle (0.1pt) node[below]{$v_5$};

            \draw[-, line width=0.4mm, lightgray, name path=l4] (A2)++(-120:0.5cm)-- ++(60:4cm) node[right, xshift=-2pt]{$\ell_4$};
        
        \end{tikzpicture}
        \caption{Minimal length $\loweru_3$ of $z_3$ in $\pentago_{\,4}$.}
        \label{fig:pent4}
    \end{subfigure}%
    \begin{subfigure}{0.4\textwidth}
    \centering
        \begin{tikzpicture}
    		\draw [dashed, name path=horizontal] (-0.5,0)--(4,0) node[at end, below, xshift=-0.75cm]{\scriptsize horizontal};
            
            \coordinate (v1) at (1,0);
            \draw[-, name path=one2] (v1)--++(120:1.5cm) node[midway, left]{$z_1$};
            \draw[-, name path=two2] (v1)++(120:1.5cm)-- ++(60:2cm) node[midway, left]{$z_2$};
            \draw[-, name path=three2] (v1)++(120:1.5cm)++(60:2cm)-- ++(-60:2.5cm);
            \draw[-, name path=four2] (v1)++(120:1.5cm)++(60:2cm)++(-60:2.5cm)-- ++(-120:1cm);
            \draw[-, name path=five2] (v1)++(120:1.5cm)++(60:2cm)++(-60:2.5cm)++(-120:1cm)-- (v1);
            
            \path [name intersections={of=one2 and two2, by={B2}}];
            \path [name intersections={of=two2 and three2, by={C2}}];
            \path [name intersections={of=three2 and four2, by={D2}}];
            \path [name intersections={of=four2 and five2, by={E2}}];
            \path [name intersections={of=five2 and one2, by={A2}}];

            \filldraw[black] (A2) circle (0.1pt) node[below, xshift=0.25cm]{$v_1$};
            \filldraw[black] (B2) circle (0.1pt) node[left]{$v_2$};
            \filldraw[black] (C2) circle (0.1pt) node[above]{$v_3$};
            \filldraw[black] (D2) circle (0.1pt) node[right]{$v_4$};
            \filldraw[black] (E2) circle (0.1pt) node[below]{$v_5$};

            % dashed lines and supporting points 
            \draw[-, line width=0.4mm, lightgray, name path=l4] (A2)++(-120:0.5cm)-- ++(60:3.5cm) node[right, xshift=-2pt]{$\ell_4$};

            \draw[dashed, name path=d1] (D2)-- ++(-60:1cm);
        
        \end{tikzpicture}
        \caption{Minimal length $\loweru_3$ of $z_3$ in $\pentago_{\,5}$.}
        \label{fig:pent5}
     \end{subfigure}%
        \caption{The pentagon special cases $4$ and $5$.}
        \label{fig:PentCases5,6}
\end{figure}
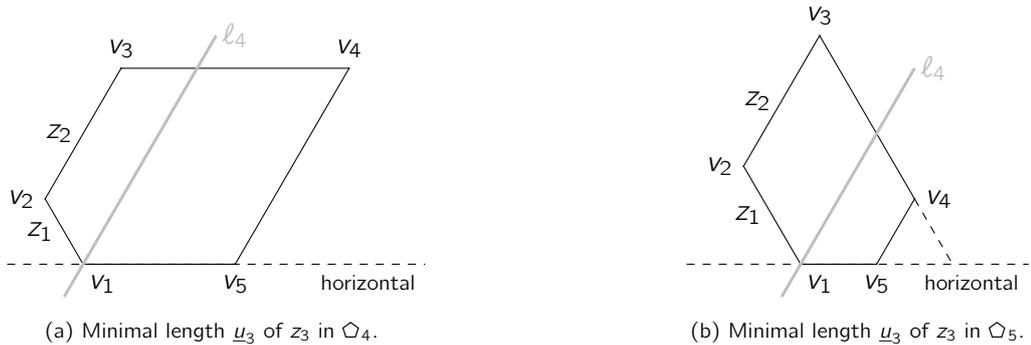
There is nothing to do in case $\pentago_{\,2.1}$, as the situation described in Section~\ref{sec:quadrilaterals} does not appear here. But it does for $\pentago_{\,2.2}$ if $z_1>z_2$. \autoref{fig:pent2.2} shows the smallest distance possible for $z_3$, that is $\ell_4$ intersecting the horizontal precisely at $v_1$.  The dashed lines therein generate two triangles, one having the vertices $v_2$, $v_3$ and the intersection of the first side of length $z_1$ with the dashed line, and the other with vertices $v_1$, $w$ and the intersection of $\ell_4$ with $\overline{v_3w}$. Both of these triangles are regular with side length $z_2$. Since the segment $\overline{v_3w}$ has length $z_1$, this yields a minimal length of $\loweru_3=z_1-z_2$ for the third side of the pentagon $\pentago_{\,2.2}$. Together with $\upperu_3=z_1$ from \eqref{eq:u3_pentagons} this yields $z_3\in(z_1-z_2,z_1)$, see Table~\ref{tab:angles_for_pentagons}.

Dealing with $\pentago_{\,4}$ and $\pentago_{\,5}$ draws analogies to the trapezoid case $\squad_{\,4}$. As illustrated in \autoref{fig:pent4}, the line $\ell_4$ (again presented such that it intersects the horizontal precisely at $v_1$ and that $z_3$ has minimal possible length) divides the pentagon into an isosceles trapezoid and a parallelogram. Due to the geometry of the trapezoid, $z_3$ must therefore be at least of length $\loweru_3=z_1$.

Turning focus to $\pentago_{\,5}$, the pentagon is divided in the same manner, see \autoref{fig:pent5}. Here, with analogous arguments, $z_3$ has to be at least of length $\loweru_3=z_1$. Together with $\upperu_3$ from \eqref{eq:u3_pentagons} this yields $z_3\in(z_1,z_1+z_2)$. 

\begin{figure}[t!]
    \centering
    \begin{subfigure}{0.5\textwidth}
        \includegraphics[scale=0.6]{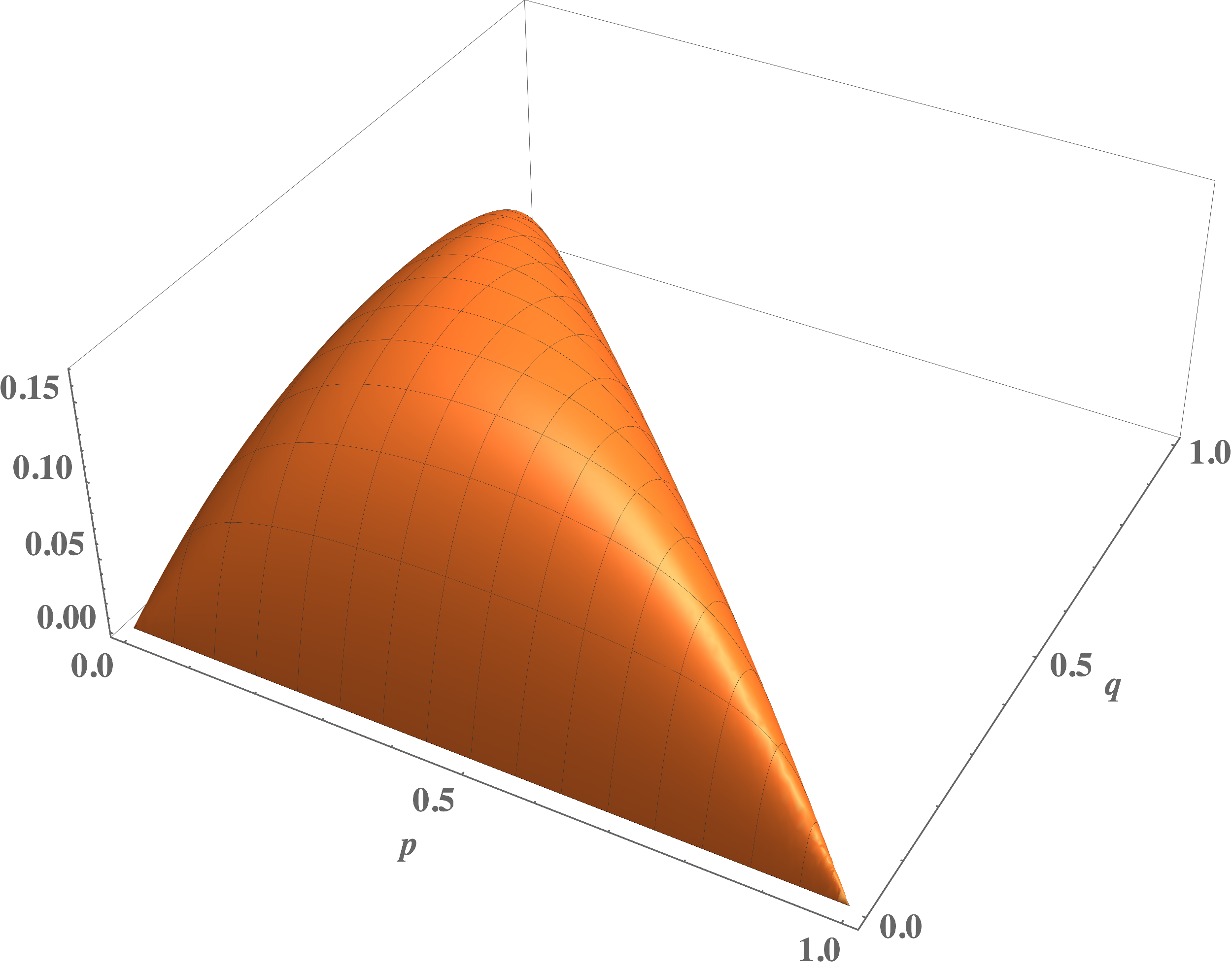}
        \caption{Plot of $\PP(N_{p,q}=5)$}
        \label{fig:plot_of_p5}
    \end{subfigure}%
    \begin{subfigure}{0.5\textwidth}
        \includegraphics[scale=0.6]{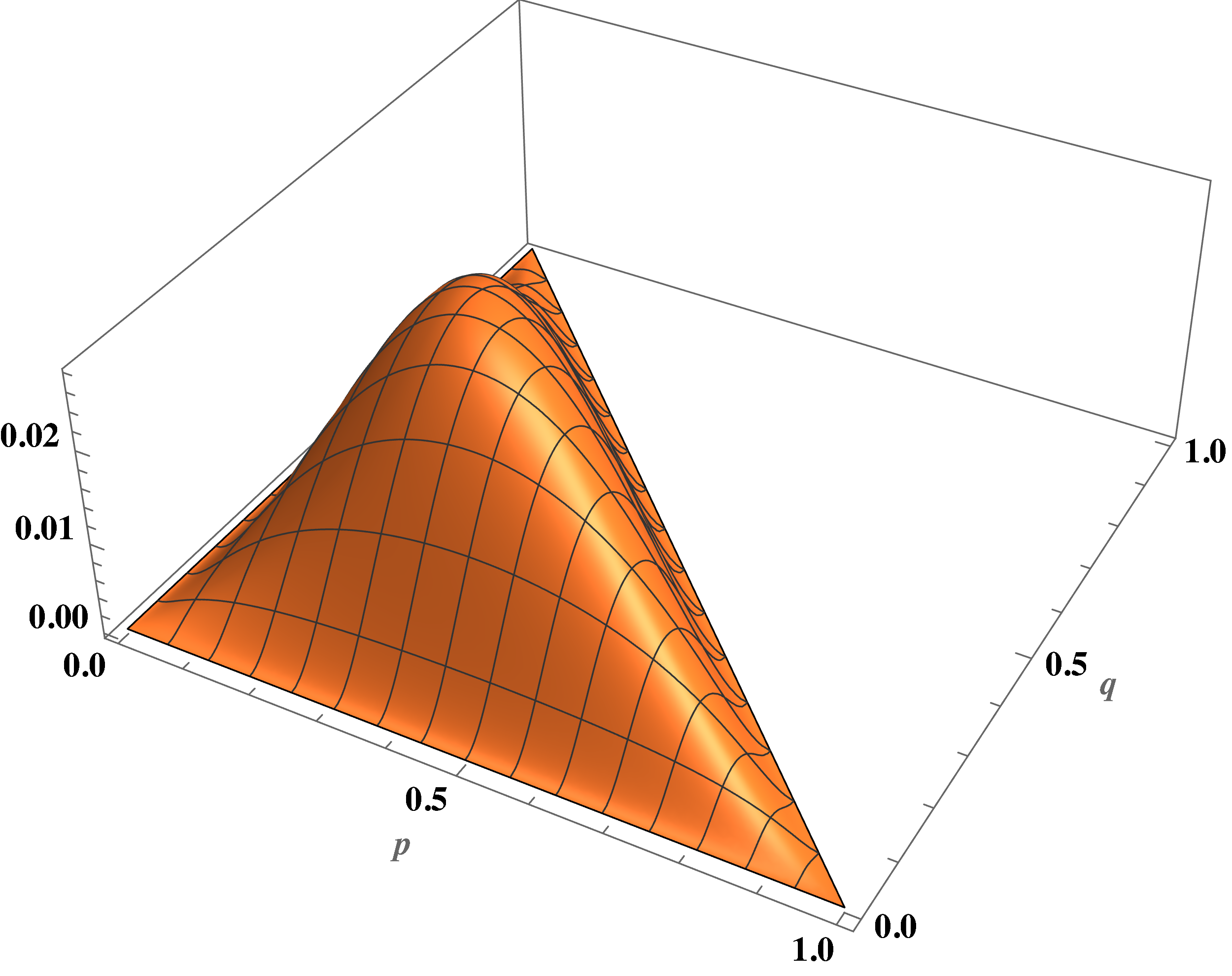}
        \caption{Plot of $\PP(N_{p,q}=6)$.}
        \label{fig:plot_of_p6}
    \end{subfigure}%
    \caption{Plots of $\PP(N_{p,q}=5)$ and $\PP(N_{p,q}=6)$.}
\end{figure}

Putting together these observations, the density \eqref{eq:density_z5} can now be integrated. For example, consider the case $\pentago_{\,5}$. We obtain
\begin{align*}
    &\PP(\tc_{p,q}=\pentago_{\,5})\\
    &=\int_{[0,\pi)^5} \int_0^\infty \int_0^\infty \int_{z_1}^{z_1+z_2} \frac{9}{8\lambda}\, \exp\bigg(
    -\frac{\sqrt{3}}{2} \big(p z_1+q( z_3-z_2) + z_2\big)\bigg)\one\big\{\varphi_0=0, \varphi_1=2\pi/3, \ldots \\[5pt]
    &\qquad \ldots\varphi_2=\pi/3,\varphi=-\pi/3,\varphi_4=-2\pi/3 \big\} \, \dint z_3\dint z_2\dint z_1G^{\otimes 5}(\dint(\varphi_0,\dots,\varphi_4)) \\[10pt]
    &={9pq^2(1-p-q)\over 8\lambda}\int_0^\infty \int_0^\infty \int_{z_1}^{z_1+z_2}\exp\bigg(
    -\frac{\sqrt{3}}{2} \big(p z_1+q( z_3-z_2) + z_2\big)\bigg)\,\dint z_3\dint z_2\dint z_1\\
    &={9pq^2(1-p-q)\over 8\lambda}{2\over\sqrt{3}\,q}\int_0^\infty \int_0^\infty \exp\Big(-{\sqrt{3}\over 2}((p+q)z_1+z_2)\Big)\Big(\exp\Big(-{\sqrt{3}\over 2}qz_2\Big)-1\Big)\,\dint z_2\dint z_1\\
    &={9pq^2(1-p-q)\over 8\lambda}{4\over 3(1-q)}\int_0^\infty \exp\Big(-{\sqrt{3}\over 2}(p+q)z_1\Big)\,\dint z_1\\
    &={9pq^2(1-p-q)\over 8\lambda}{8\over 3\sqrt{3}(1-q)(p+q)}\\
    & = \frac{3 p q^2 (-p-q+1)^2}{(3-3 q) (p+q) \left(-p^2-p q+p-q^2+q\right)},
\end{align*}
where the calculations where carried out similar to Sections~\ref{sec:triangles} and~\ref{sec:quadrilaterals}. Dealing with the remaining cases in the same way eventually leads to the following result, see also \autoref{fig:plot_of_p5}.

\begin{lem}
    \label{lem:proportion_pentagons}
    In the setup of Theorem \autoref{thm:main_theorem}, for all $0<p,q<1$ with $0<p+q<1$, we have that
    \begin{align*}
        \beta_{p,q}^{-1} \, \big[ \, 2pq(1-p-q)(p^2+q^2) &- 6p^2q^2(p+q)(1-p-q) \\[5pt]
            &- 2pq(p+q)(1-p-q)  + 8p^2q^2(1-p-q) \, \big].
    \end{align*}
The maximal value for $\PP(N_{p,q}=5)$ is attained precisely if $p=q=1/3$ and is given by
    \begin{equation*}
    	\max_{0<p+q<1} \PP(N_{p,q}=5)=\PP(N_{1/3,1/3}=5)=\frac16.
    \end{equation*}
\end{lem}

% ====================================================================
\subsection{The hexagon case}
\label{sec:hexagons}
% ====================================================================

We finally deal with the probability that the typical cell is a hexagon. One can easily see that there is only one possible combination of angles that can lead to such a shape, see Table~\ref{tab:angles_for_hexagons}.

\begin{figure}[b!]
    \centering
     \begin{subfigure}{0.45\textwidth}
        \centering
        \begin{tikzpicture}
    		\draw [dashed, name path=horizontal] (-1,0)--(4.5,0) node[at end, below, xshift=-0.75cm]{\scriptsize horizontal};
            
            \coordinate (v1) at (0,0);
            \draw[-, name path=one2] (v1)--++(120:2.5cm) node[midway, left]{$z_1$};
            \draw[-, name path=two2] (v1)++(120:2.5cm)-- ++(60:2cm) node[midway, left]{$z_2$};
            \draw[-, name path=three2] (v1)++(120:2.5cm)++(60:2cm)-- ++(0:1.5cm);
            \draw[-, name path=four2] (v1)++(120:2.5cm)++(60:2cm)++(0:1.5cm)-- ++(-60:2.5cm);
            \draw[-, name path=five2] (v1)++(120:2.5cm)++(60:2cm)++(0:1.5cm)++(-60:2.5cm)-- ++(-120:2cm);
            \draw[-, name path=six2] (v1)++(120:2.5cm)++(60:2cm)++(0:1.5cm)++(-60:2.5cm)++(-120:2cm)--(v1);

            \path [name intersections={of=one2 and two2, by={B2}}];
            \path [name intersections={of=two2 and three2, by={C2}}];
            \path [name intersections={of=three2 and four2, by={D2}}];
            \path [name intersections={of=four2 and five2, by={E2}}];
            \path [name intersections={of=five2 and six2, by={F2}}];
            \path [name intersections={of=six2 and one2, by={A2}}];

            \filldraw[black] (A2) circle (0.1pt) node[below, xshift=0.25cm]{$v_1$};
            \filldraw[black] (B2) circle (0.1pt) node[left]{$v_2$};
            \filldraw[black] (C2) circle (0.1pt) node[left]{$v_3$};
            \filldraw[black] (D2) circle (0.1pt) node[right]{$v_4$};
            \filldraw[black] (E2) circle (0.1pt) node[right]{$v_5$};
            \filldraw[black] (F2) circle (0.1pt) node[below]{$v_6$};

            \draw[-, line width=0.4mm, lightgray, name path=l4] (A2)++(-120:0.5cm)-- ++(60:5cm) node[right, xshift=-2pt]{$\ell_5$};

            % dashed lines and supporting points
            \draw[dashed] (C2)-- ++(60:1.5cm);
            \draw[dashed] (C2)++(60:1.5cm)--++(-60:1.5cm);
            \draw[dashed] (E2)--++(-60:2cm);
            
        \end{tikzpicture}
        \caption{Minimal length $\loweru_4$ of $z_4$ in $\hexago_{\,1.1}$.}
        \label{fig:hex1}
    \end{subfigure}%
    \begin{subfigure}{0.55\textwidth}
        \centering
        \begin{tikzpicture}
    		\draw [dashed, name path=horizontal] (-1,0)--(6,0) node[at end, below, xshift=-0.75cm]{\scriptsize horizontal};
            
            \coordinate (v1) at (0,0);
            \draw[-, name path=one2] (v1)--++(120:2.5cm) node[midway, left]{$z_1$};
            \draw[-, name path=two2] (v1)++(120:2.5cm)-- ++(60:2cm) node[midway, left]{$z_2$};
            \draw[-, name path=three2] (v1)++(120:2.5cm)++(60:2cm)-- ++(0:3.5cm);
            \draw[-, name path=four2] (v1)++(120:2.5cm)++(60:2cm)++(0:3.5cm)-- ++(-60:2.5cm);
            \draw[-, name path=five2] (v1)++(120:2.5cm)++(60:2cm)++(0:3.5cm)++(-60:2.5cm)-- ++(-120:2cm);
            \draw[-, name path=six2] (v1)++(120:2.5cm)++(60:2cm)++(0:3.5cm)++(-60:2.5cm)++(-120:2cm)--(v1);

            \path [name intersections={of=one2 and two2, by={B2}}];
            \path [name intersections={of=two2 and three2, by={C2}}];
            \path [name intersections={of=three2 and four2, by={D2}}];
            \path [name intersections={of=four2 and five2, by={E2}}];
            \path [name intersections={of=five2 and six2, by={F2}}];
            \path [name intersections={of=six2 and one2, by={A2}}];

            \filldraw[black] (A2) circle (0.1pt) node[below, xshift=0.25cm]{$v_1$};
            \filldraw[black] (B2) circle (0.1pt) node[left]{$v_2$};
            \filldraw[black] (C2) circle (0.1pt) node[left]{$v_3$};
            \filldraw[black] (D2) circle (0.1pt) node[right]{$v_4$};
            \filldraw[black] (E2) circle (0.1pt) node[right]{$v_5$};
            \filldraw[black] (F2) circle (0.1pt) node[below]{$v_6$};

            \draw[-, line width=0.4mm, lightgray, name path=l4] (A2)++(-120:0.5cm)-- ++(60:6cm) node[right, xshift=-2pt]{$\ell_5$};

            % dashed lines and supporting points
            \draw[dashed] (E2)--++(-60:2cm);
        \end{tikzpicture}
        \caption{Minimal length $\loweru_4$ of $z_4$ in $\hexago_{\,1.2}$.}
        \label{fig:hex2}
     \end{subfigure}%
        \caption{The two hexagon cases.}
        \label{fig:hex}
\end{figure}
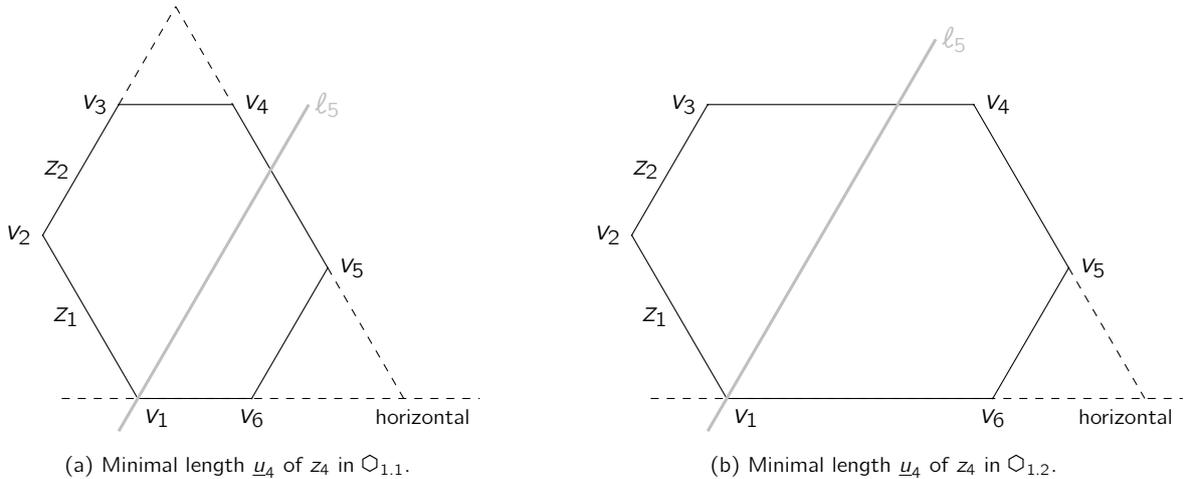

Using \eqref{eq:density_z^n} with $n=6$, it follows that the joint density of $(\Phi_0,\dots,\Phi_5,Z_1,\dots,Z_4)$ is given by
\begin{align}
\label{eq:density_z6}
    (\varphi_0,\dots,\varphi_5,z_1,\dots,z_4) &\longmapsto
    \frac{9\sqrt{3}}{16\lambda}\, \exp\bigg(-\frac{1}{2} \sum_{i=1}^6 z_i\lambda(\varphi_i)\bigg) \nonumber \\[5pt]
    &\hspace{2cm} \times \one\big\{ P(\varphi_0,\varphi_1,\varphi_2,\varphi_3,\varphi_4,\varphi_5,z_1,z_2,z_3,z_4)\in\mathsf{poly}_6\big\}.
\end{align}
The upper integral limits for the side lengths $z_i$ are again given by \eqref{eq:u_i} and equal $\upperu_1,\upperu_2,\upperu_3=\infty$ and $\upperu_4=z_1+z_2$. Similar to the prior sections, we have to be careful with the lower integration limits for some of the side lengths. Here, we have to ensure that $\ell_5$, the line corresponding to the hexagon side $z_5$ does not intersect the first side of length $z_1$ prior to intersecting the horizontal, see \autoref{fig:hex}. For these lower integration limits $\loweru_i$, adopting our prior notation, we subdivide the hexagon case into two subcases denoted by $\hexago_{\,1.1}$ and $\hexago_{\,1.2}$, respectively. In the first situation, we restrict $z_3<z_1$, in the latter we let $z_3>z_1$, see Figures~\ref{fig:hex1} and~\ref{fig:hex2} for an illustration. For $\hexago_{\,1.1}$, this leads to $\loweru_4=z_1-z_3$, and for $\hexago_{\,1.2}$ we can allow $\loweru_4=0$. This can be clarified by consulting \autoref{fig:hex} with similar geometric arguments as in the prior sections. Note, that in \autoref{fig:hex2} the position of $\ell_5$, which is again chosen in the minimal way such that it intersects the horizontal precisely at $v_1$, indicates that $v_4$ can be arbitrarily small due to the already ensured length of $z_3$ to be larger than $z_1$. 
\begin{table}[t!]
    \centering
     \renewcommand{\arraystretch}{1.5} 
     \begin{tabular}{p{0.65cm}|c|c|c|c|c|c|c|c|c|c}
     & $\varphi_0$ & $\varphi_1$ & $z_1$ & $\varphi_2$ & $z_2$ & $\varphi_3$ & $z_3$ & $\varphi_4$ & $z_4$ & $\varphi_5$ \\
     \hline
     \begin{tabular}{@{}c@{}}
            $\hexago_{\,1.1}$ \\
            $\hexago_{\,1.2}$ \\
        \end{tabular}
        & $0$ & {\Large $\frac{2\pi}{3}$} &
        \begin{tabular}{@{}c@{}}
            $(0,\infty)$ \\
            $(0,\infty)$ \\
        \end{tabular}
        & {\Large $\frac\pi3$} &
        \begin{tabular}{@{}c@{}}
            $(0,\infty)$ \\
            $(0,\infty)$ \\
        \end{tabular}
        & $0$ & 
        \begin{tabular}{@{}c@{}}
            $(0,z_1)$ \\
            $(z_1,\infty)$ \\
        \end{tabular}
        & {\Large $\frac{-\pi}{3}$} &
        \begin{tabular}{@{}c@{}}
            $(z_1-z_3, z_1+z_2)$ \\
            $(0,z_1+z_2)$ \\
        \end{tabular}
        & {\Large $\frac{-2\pi}{3}$ } \\
    \end{tabular}
    \renewcommand{\arraystretch}{1} %back to normal 
    \caption{Values of angles $\varphi_i$ and intervals $(\loweru_i,\upperu_i)$ for side lengths $z_i$ resulting in hexagons.}
    \label{tab:angles_for_hexagons}
\end{table}

This can now be used to integrate the density in \eqref{eq:density_z6}, which eventually leads to the following result.

\begin{lem}
    \label{lem:proportion_hexagons}
    In the setup of Theorem \autoref{thm:main_theorem}, for all $0<p,q<1$ with $0<p+q<1$, we have that
    \begin{align*}
        \beta_{p,q}^{-1} \, \big[ \, -2 p^2 q^2 (1-p-q)^2 \, \big].
    \end{align*}
The maximal value for $\PP(N_{p,q}=6)$ is attained precisely if $p=q=1/3$ and is given by
\begin{equation*}
	\max_{0<p+q<1} \PP(N_{p,q}=6)=\PP(N_{1/3,1/3}=6)=\frac{1}{36}.
\end{equation*}
\end{lem}

\pagebreak 

\subsection*{Acknowledgement}
CT has been supportet by the DFG priority program SPP 2265 \textit{Random Geometric Systems}. 
We are grateful to Tom Kaufmann for inspiring ideas and constructive discussions on the subject of this paper.

% ====================================================================
\printbibliography
% ====================================================================

\end{document}